\documentclass[a4paper,10pt]{article}
\usepackage[utf8x]{inputenc}
\usepackage{amsmath}
\usepackage{amsfonts}
\usepackage{amssymb}
\usepackage{amsthm}
\date{} 
\usepackage{epsfig}
\usepackage{mathrsfs}
\usepackage{color}
\usepackage[british,english]{babel}
\usepackage{amsthm}
\newtheorem{prop}{\quad Proposition}[section]
\newtheorem{theoreme}{\quad Theorem}[section]
\newtheorem{lem}{\quad Lemma}[section]
\newtheorem{rem}{\quad Remark}[section]
\newtheorem{cor}{\quad Corollary}
\def \1{\textrm{\dsrom{1}}}

\begin{document}
%\newtheorem{prop}{Proposition}
%\newtheorem{lem}{Lemma}

% the locale $L^2$ norm, of velocity field, solution of Navier equation
\title{ Insensitizing controls for the Navier-Stokes equations}
\author{Mamadou Gueye  \thanks{ Universit\'e
Pierre et Marie Curie-Paris 6, UMR 7598 Laboratoire Jacques-Louis Lions, bo\^ite courrier 187, 4 place Jussieu 75252 Paris Cedex 05, France. E-mail: {\bf gueye@ann.jussieu.fr}.}}

\maketitle

\begin{abstract}
In this paper, we deal with the existence of insensitizing controls for the Navier-Stokes equations in a bounded domain with Dirichlet boundary conditions. We prove that there exist controls insensitizing the $L^2$ -norm of the observation of the solution in an open subset $\mathcal{O}$
of the domain, under suitable assumptions on the data.
 This problem is equivalent to an exact controllability result for a cascade system. First we prove a global Carleman inequality for the linearized Navier-Stokes system with right-hand side, which leads to the null controllability at any time $T>0$. Then, we deduce a local null controllability result for the cascade system.  
\end{abstract}

\section{Introduction}
Let $\Omega\subset \mathbb R^N$($N=2\ \text{or}\ 3$) be a bounded connected open set whose boundary $\partial \Omega$ is regular enough (for instance of class $C^2$). Let $\omega$ and $\mathcal O$ be two open and nonempty subsets of $\Omega$ (resp. the control domain and the observatory)  and let $T>0$. We will use the notation $Q=\Omega\times (0,T)$ and $\Sigma=\partial\Omega\times (0,T)$. $C$ stands for a generic constant which depends only on $\Omega$, $\omega$, $\mathcal O$ and $T$. \\
%$\omega\cap \mathcal O \neq \emptyset$, 
The Navier-Stokes equations describe the motion of an incompressible fluid such as water, air, oil... They appear in the study of many phenomena, either alone or coupled with other equations. For instance, they are used in theoretical studies in meteorology, in aeronautical sciences, in environmental sciences, in plasma physics, in the petroleum industry, etc.\\%%These equations are to be solved for an unknown velocity vector $y(x,t)=(y_i(x,t)_{1\leq i\leq N})\in \mathbb R^N$ and pressure $p(x,t)\in \mathbb R$, defined for position $x\in \Omega$ and time $t\in [0,T]$. 
%Here, we consider this system with distributed controls and incomplete data
First let us recall some usual spaces in the context of Navier-Stokes equations:
$$V=\left\{y\in H^1_0(\Omega)^N \ ;\nabla \cdot y=0\ \text{in}\ \Omega\right\},$$
and
$$H=\left\{y\in L^2(\Omega)^N \ ;\nabla \cdot y=0\ \text{in}\ \Omega, \ y\cdot n=0\ \text{on}\ \partial \Omega\right\}.$$

To be more specific about the investigated problem, we introduce the following control system with incomplete data
\begin{equation} \label{1} \left\lbrace { \begin{tabular}{llll}
$y_t- \Delta y + (y,\nabla)y+\nabla p=f+v\mathbf{1}_{\omega}$ & \textrm{in }$Q$,\\
$ \nabla \cdot y=0 $ & \textrm{in }$Q$,\\
$y=0$ & \textrm{on }$\Sigma$,\\
$  y_{|t=0}=y^0+\tau \widehat y^0 $ & \textrm{in }$\Omega$.\\ \end{tabular}
} \right. 
\end{equation}
% \begin{equation} \label{1_jp} \left\lbrace { \begin{array}{lll}
%  y_t- \Delta y + (y,\nabla)y+\nabla p=f+v\mathbf{ 1}_{\omega},\ \nabla \cdot y=0  & \textrm{in }Q,\\
% y=0 & \textrm{on }\Sigma,\\
%   y_{|t=0}=y^0+\tau \widehat y^0  & \textrm{in } \Omega.\\ \end{tabular}
% } \right. 
% \end{equation}

Here, $y(x,t)=(y_i(x,t))_{1\leq i\leq N}$ is the velocity of the particles of an incompressible fluid,
$v$ is a distributed control localized in $\omega$, $f(x,t)=(f_i(x,t))_{1\leq i\leq N}\in L^2(Q)^N$ is a given, externally applied force, and the initial state  $y_{|t=0}$ is partially unknown. We suppose that $y^0\in H$, $\widehat y^0\in H$ is unknown with $\|\widehat y^0\|_{L^2(\Omega)^N}=1$ and that $\tau$ is a small unknown real number.\\
The aim of this paper is to prove the existence of controls that insensitize some functional $J_{\tau}$ (the sentinel) depending on the velocity field $y$. That is to say, we have to find a control $v$ such that the influence of the unknown data $\tau\widehat y^0$ is not perceptible for our sentinel:  \\
\begin{equation}\label{3}
\left. \dfrac{\partial J_{\tau}(y)}{\partial \tau}\right|_{\tau=0}=0\quad \forall \widehat y^0\in L^2(\Omega)^N \quad \text{such that}\quad \|\widehat y^0\|_{L^2(\Omega)^N}=1.
\end{equation}
In the pioneering work \cite{llions}, J.-L. Lions considers this kind of problem and introduces many related questions. One of these questions, in non-classical terms, was the existence of insensitizing controls for the Navier-Stokes equations (see \cite{llions}, page 56).\\
In the literature the usual sentinel is given by the square of the local $L^2$-norm of the state variable $y$, on which we will be interested here:
\begin{equation}\label{4}
 J_{\tau}(y)={1\over2}\iint\limits_{\mathcal O\times (0,T)}|y|^2\textrm{d}x\textrm{d}t.
\end{equation}
However, in \cite{guerrero}, the author considers the gradient of the state for a linear heat system with potentials and more recently in \cite{Gue} the same author treats the case of the curl of the solution for a Stokes system. Here we will focus on the nonlinear control problem of insensitizing the Navier-Stokes equations. \\

%It is by now classical (see \cite{Temam}), since J. Leray's work in the '30s, that system (\ref{1}) admits a unique weak solution under a smallness assumption on the initial data. Futhermore, provided that the data is sufficiently smooth (see \cite{Temam}, for instance, §3), the initial boundary-value problem (\ref{1}) for the Naviers-Stokes equations possesses a unique smooth solution on some interval of time $(0,T_{\star})$. 
%%See also generic solvability of the Naviers-Stokes by A. V. Fursikov \cite{Fu} for a partial answer to the question to know whether system (\ref{1}) is solvable for an arbitrary pair $y_{|t=0}$, $f$.\\

The special form of the sentinel $ J_{\tau}$ allows us to reformulate our insensitizing problem as a controllability problem of a cascade system (for more details, see \cite{Bodart}, for instance). In particular, condition (\ref{3}) is equivalent to $z_{|t=0}=0$ in $\Omega$, where $z$ together with $w$ solves the following coupled system:
\begin{equation} \label{2} \left\lbrace { \begin{tabular}{llll}
 $w_t- \Delta w + (w,\nabla)w+\nabla p^0=f+v\mathbf{ 1}_{\omega},\ \nabla \cdot w=0 $ & \textrm{in }$Q$,\\
$-z_t- \Delta z + (z,\nabla^t)w- (w,\nabla)z+\nabla q=w\mathbf{ 1}_{\mathcal O},\ \nabla \cdot z=0 $ & \textrm{in }$Q$,\\
$w=z=0$ & \textrm{on }$\Sigma$,\\
$  w_{|t=0}=y^0,\  z_{|t=T}=0 $ & \textrm{in }$\Omega$.\\ \end{tabular}
} \right. 
\end{equation}
Here, $(w,p^0)$ is the solution of system (\ref{1}) for $\tau=0$, the equation of $z$ corresponds to a formal adjoint of the equation satisfied by the derivative of $y$ with respect to $\tau$ at $\tau=0$ (see \eqref{tau} below) and we have denoted
$$((z,\nabla^t)w)_i= \sum_{j=1}^Nz_j\partial_iw_j\quad i=1,\dotsc,N.$$
Indeed, differentiating $y$ solution of \eqref{1} with respect to $\tau$ and evaluating it at $\tau=0$, condition \eqref{3} reads
\begin{equation}
 \iint\limits_{\mathcal O\times (0,T)}wy_\tau\textrm{d}x\textrm{d}t=0\quad\forall \widehat y^0\in L^2(\Omega)^N \quad \text{such that}\quad \|\widehat y^0\|_{L^2(\Omega)^N}=1,
\end{equation}
where $y_\tau$ is the derivative of $y$ solution of \eqref{1} at $\tau=0$.
Then, $y_{\tau}$ solves 
\begin{equation} \label{tau} \left\lbrace { \begin{tabular}{llll}
$y_{\tau,t}- \Delta y_\tau + (y_\tau,\nabla)y+(y,\nabla)y_{\tau}+\nabla p_\tau=0$ & \textrm{in }$Q$,\\
$ \nabla \cdot y_\tau=0 $ & \textrm{in }$Q$,\\
$y_\tau=0$ & \textrm{on }$\Sigma$,\\
$  y_{\tau|t=0}= \widehat y^0 $ & \textrm{in }$\Omega$.\\ \end{tabular}
} \right. 
\end{equation}
Hence, substituting $w\mathbf{ 1}_{\mathcal O}$ by the left-hand side of the equation of $z$ in \eqref{2} and integrating by parts we obtain
\begin{equation}
 \int\limits\limits_{\Omega}z_{|t=0}\widehat y^0\textrm{d}x=\iint\limits_{\mathcal O\times (0,T)}wy_\tau\textrm{d}x\textrm{d}t\quad \forall \widehat y^0\in L^2(\Omega)^N \quad \text{such that}\quad \|\widehat y^0\|_{L^2(\Omega)^N}=1.
\end{equation}

We will prove the following controllability result for system (\ref{2}):
\begin{theoreme}
Let $m>5$ be a real number and $y^0=0$. Assume that $\omega\cap \mathcal O \neq \emptyset$. Then, there exist $\delta>0$ and $C^{\star}>0$ depending on $\omega$, $\Omega$, $\mathcal O$ and $T$ such that for any $f\in L^2(Q)^N$ satisfying $\|e^{C^{\star}/t^m}f\|_{ L^2(Q)^N}<\delta$, there exists a control $v\in L^2(\omega\times(0 ,T))^N$ and a corresponding solution $(w,z)$ to (\ref{2}) satisfying  $z_{|t=0}=0$ in $\Omega$. 
\end{theoreme}
\begin{rem}
 Furthermore, in addition to insensitizing the functional $J_{\tau}$ one can steer the state $w$ to $0$ at time $t=T$  just by paying an extra condition on $f$ at time $t=T$
\begin{equation}
  \|e^{C_{\star}/t^m(T-t)^m}f\|_{ L^2(Q)^N}< +\infty,
\end{equation}
for  a  constant $C_{\star}$ that maybe different to the one given in Theorem $1.1$.
\end{rem}
\begin{rem}
The condition $y_0=0$ in the main theorem is due to the fact that the first equation in (\ref{2}) is forward and the second one is backward in time. Most of the insensitizing works in the parabolic case, even for linear equations, assume this condition on the initial data. A study of the possible initial conditions which can be insensitized is made
for the heat equation in \cite{ter3}.
% $y_0\neq0$ constituting an obstruction to the null controllability at time $t=0$ for the solution of the second equation in (\ref{2}). Notice also that this is consistant with the requirement on the source term. 
 %%It seems natural that insensitizing (\ref{4}) with respect to small disturbances near the %%equilibrium is the relevant property; this explains the condition $y_0=0$ in the main theorem.
% Consequently, we do not know whether Theorem $1.1$ still holds if $y_0$ is not identically zero.
This work suggests that the answer is not obvious.
\end{rem}

As announced, we have the following result.
\begin{cor}
 There exists insensitizing controls $v$ for the functional $J_{\tau}$ given by (\ref{4}). 
\end{cor}
%%The proof of this result is by now a classical consequence of Theorem $1.1$ (see, for instance, \cite{Bodart}).\\\\
Before going further, let us recall some of the results avalaible in the literature. 
Most known results concerning insensitizing controls are for parabolic systems. Nevertheless, one can cite the results in \cite{dag}  for the 1-D wave equation. In \cite{Tebou}, the controllability of coupled wave equations is studied.\\
%%In the context of nonlinear problems it seems natural to relax condition (\ref{4}) since exact controllability is very hard to achieve
In order to get rid of the condition $y^0=0$, in \cite{Bodart}, the authors consider $\epsilon$-insensitizing controls (i.e., $v$ such that $|\partial_\tau J_\tau(y)_{|\tau=0}|\leq\epsilon$ for all $\epsilon>0$ ) for the semilinear heat system, with $C^1$ and globally Lipschitz nonlinearities, and prove that this condition is equivalent to an approximate controllability result for a cascade system which is established therein. In [9], condition $y^0=0$ has been removed for the linear heat equation
when ${\mathcal O}\subset\omega$ and when ${\mathcal O}=\Omega$, if this is not the case, some negative results are also provided. In \cite{ter1}, the author proves the existence of insensitizing controls for the same semilinear heat system. This last result is extended in \cite{Gonzales} to super-linear nonlinearities.\\
For parabolic systems arising from fluids dynamics the first attempt to treat the insensitizing problem is \cite{Car} for a large scale ocean circulation model (linear). In \cite{Gue}, as we have already mentioned, the author treats both the case of a sentinel given by $L^2$-norm of the state and $L^2$-norm of the curl of the state of a linear Stokes system. 

As long as insensitizing controls have been considered the condition $\omega\cap \mathcal O \neq \emptyset$ has always been imposed. But, from \cite{ter2} and \cite{dter}, we see that this is not a necessary condition for $\epsilon$-insensitizing controls. For instance, the authors have proved in \cite{dter} that there exists $\epsilon$-insensitizing controls of $J_{\tau}$ for linear heat equations with no intersecting observation and control regions in one space dimension using the spectral theory.\\
Furthermore, the insensitizing problem, as we have seen in this special case, is directly related to control problems for coupled systems. In particular, one could ask whether it is possible to control both states of a coupled system just by acting on one equation. In this spirit, the authors in \cite{CoronGue} show some controllability results for the Navier-Stokes equations with controls having a vanishing component.
In \cite{Gue} and \cite{guerrero}, as well as some insensitizing problems, the author studied this problem respectively for Stokes and heat systems in a more general framework. Also, for more general coupled parabolic systems with only one control force, some results are available in \cite{Gonzalez} and \cite{Ammar}.\\
Finally, recently in \cite{Yan} the existence of insensitizing controls for a forward stochastic heat equation was proved by means of some global Carleman estimates.\\
The rest of the paper is organized as follows. In Section 2 we give some results which will be useful for our purpose. In Section 3 we prove the Carleman estimate. In Section 4 we treat the linear case. Finally, in Section 5 we prove the insensitivity result.
\section{Technical results}
In the context of the null controllability analysis of parabolic systems, Carleman estimates are a very powerful tool (see \cite{FI},\cite{immanuvilov},\cite{Fer},...).
In order to state our Carleman estimate, we need to define some weight functions. 
Let $\omega_0$ be a nonempty open subset of $\omega\cap\mathcal O$, and set :

\begin{equation}\label{weight}
{\alpha_m(x,t)={\exp(\lambda k{m+1\over m}\|\eta^0\|_{\infty})-\exp{\lambda(k\|\eta^0\|_{\infty}+\eta^0(x))}\over t^m(T-t)^m},\quad\xi_m(x,t)={\exp{\lambda(k\|\eta^0\|_{\infty}+\eta^0(x))}\over t^m(T-t)^m}},
 \end{equation}
%and  we shall use the notation
%\begin{equation}
%\alpha_m^{\star}(t)=\max_{x\in \overline \Omega}\alpha_m(x,t),\quad  \widehat\alpha_m(t)=\min_{x\in \overline \Omega}\alpha_m(x,t),\quad\xi_m^{\star}(t)=\min_{x\in \overline \Omega}\xi_m(x,t), \quad \widehat\xi_m(t)=\max_{x\in \overline \Omega}\xi_m(x,t).
% \end{equation}
for some parameter $\lambda>0$. 
Here, $m > 4$ and $k> m$ are fixed and $\eta^0\in C^2(\overline{\Omega})$ stands for a function that satisfies :
\begin{equation}\label{omega}
 |\nabla \eta^0|\geq K>0 \quad \text{in}\ \Omega \backslash\overline{\omega}_0,\ \eta^0>0 \ \text{in}\ \Omega \quad \text{and} \quad \eta^0=0\ \text{on} \ \partial \Omega .
\end{equation}
The proof of the existence of such a function $\eta_0$ can be found in \cite{FI}.
This kind of weight functions was also used in \cite{immanuvilov}. In the sequel, for convenience, we will fix $m=5$ and $k=10$. Thus, our weight functions read

\begin{equation}\label{weight1}
{\alpha(x,t)={\exp(12\lambda \|\eta^0\|_{\infty})-\exp{\lambda(10\|\eta^0\|_{\infty}+\eta^0(x))}\over t^5(T-t)^5},\quad\xi(x,t)={\exp{\lambda(10\|\eta^0\|_{\infty}+\eta^0(x))}\over t^5(T-t)^5}}
 \end{equation}
and  we shall use the notation
\begin{equation}
\alpha^{\star}(t)=\max_{x\in \overline \Omega}\alpha(x,t),\quad  \widehat\alpha(t)=\min_{x\in \overline \Omega}\alpha(x,t),\quad\xi^{\star}(t)=\min_{x\in \overline \Omega}\xi(x,t), \quad \widehat\xi(t)=\max_{x\in \overline \Omega}\xi(x,t).
 \end{equation}
We also introduce the following quantities:
\begin{equation}\label{io}
 I_0(s,\lambda ;u)=s^3 \lambda^4\iint\limits_{Q}e^{-2s\alpha}\xi^3|u|^2\textrm{d}x\textrm{d}t+s\lambda^2\iint\limits_{Q}e^{-2s\alpha}\xi|\nabla u|^2\textrm{d}x\textrm{d}t,
\end{equation}
\begin{equation}\label{ione}
I_1(s, \lambda;u)=s^3 \lambda^4\iint\limits_{Q}e^{-5s\alpha}\xi^3|u|^2\textrm{d}x\textrm{d}t+s\lambda^2\iint\limits_{Q}e^{-5s\alpha}\xi|\nabla u|^2\textrm{d}x\textrm{d}t+s^{-1}\iint\limits_{Q}e^{-5s\alpha}(\xi)^{-1}|\Delta u|^2\textrm{d}x\textrm{d}t,
 \end{equation}
\begin{equation}\label{itilde}
\widetilde I(s,\lambda;u)=s^3 \lambda^4\iint\limits_{Q}e^{-2s\alpha}e^{-2s\alpha^{\star}}\xi^3|u|^2\textrm{d}x\textrm{d}t+
s\lambda^2\iint\limits_{Q}e^{-2s\alpha}e^{-2s\alpha^{\star}}\xi|\nabla u|^2\textrm{d}x\textrm{d}t,
\end{equation}
for some parameter $s>0$.\\

First we state a Carleman-type estimate which holds for energy solutions of heat equations with non-homogeneous Neumann boundary conditions:
\begin{lem}
 Let us assume that $u^0\in L^2(\Omega)$, $f_1\in L^2(Q)$, $f_2\in L^2(Q)^N$, $f_3\in L^2(\Sigma)$. Then, there exists a constant $C(\Omega,\omega_0)>0$ such that the (weak) solution of 

\begin{equation}\label{eq1}
   \begin{cases}
       u_t - \Delta u =f_1+\nabla \cdot f_2  &\quad \text{in} \  Q,\\
      \displaystyle {\partial u\over \partial n}+f_2\cdot n=f_3 &\quad \text{on}\  \Sigma, \displaystyle\\
       u_{|t=0}=u^0 &\quad \text{in} \ \Omega,
   \end{cases}
\end{equation}
satisfies
\begin{multline}
I_0(s,\lambda ;u)\leq C\left(s^3 \lambda^4\iint\limits_{\omega_0\times (0,T)}e^{-2s\alpha}\xi^3|u|^2\textrm{d}x\textrm{d}t+\iint\limits_{Q}e^{-2s\alpha}|f_1|^2\textrm{d}x\textrm{d}t     \right. \\
 \left.+s^2 \lambda^2\iint\limits_{Q}e^{-2s\alpha}\xi^2|f_2|^2\textrm{d}x\textrm{d}t+s\lambda\iint\limits_{\Sigma}e^{-2s\alpha^{\star}}\xi^{\star}|f_3|^2d\sigma \textrm{d}t \right),
\end{multline}
for any $\lambda \geq C$ and $s\geq C(T^{9}+T^{10})$.
\end{lem}
In a similar form, this lemma was proved in \cite{Gon}, but with the weight defined in (\ref{weight}) for $m=1$. In order to prove Lemma $2.1$, one can follow the steps of the proof in \cite{Gon}, just taking into account that
\begin{equation}\label{alpha}
|\alpha_t|\leq KT \xi^{6/5},\quad |\alpha^{\star}_t|\leq KT (\xi^{\star})^{6/5},\quad |\alpha_{tt}|\leq KT^2\xi^{7/5},\quad \text{and}\quad |\alpha^{\star}_{tt}|\leq KT^2(\xi^{\star})^{7/5},
 \end{equation}
for some constant $K$ independent of $s$, $\lambda$ and $T$.\\

The second estimate we give here holds for solutions of Stokes systems with homogeneous Dirichlet boundary conditions:
\begin{lem}
Let us assume that $u^0\in V$, $f_4\in L^2(Q)^N$. Then, there exists a constant $C(\Omega,\omega_0)>0$ such that the solution $(u,p)\in (L^2(0,T;H^2(\Omega)^N\cap V)\cap L^{\infty}(0,T;V)\times L^2(0,T;H^1(\Omega))$, with $\int\limits_{\omega_0}p(t,x)dx=0$, of
\begin{equation}\label{eq3}
   \begin{cases}
       u_t - \Delta u +\nabla p=f_4 ,\nabla \cdot u=0 &\quad \text{in} \  Q,\\
       u=0 &\quad \text{on}\  \Sigma, \\
       u_{|t=0}=u^0 &\quad \text{in} \ \Omega,
   \end{cases}
\end{equation}
satisfies
\begin{equation}\label{eq4}
 I_0(s,\lambda ;u)\leq C\left(s^{16} \lambda^{40}\iint\limits_{\omega_0\times (0,T)}e^{-8s\widehat\alpha+6s\alpha^{\star}}(\widehat\xi)^{16}|u|^2\textrm{d}x\textrm{d}t
+s^{15/2} \lambda^{20}\iint\limits_{Q}e^{-4s\widehat\alpha+2s\alpha^{\star}}(\widehat\xi)^{15/2}|f_4|^2\textrm{d}x\textrm{d}t\right),
\end{equation}
for any $\lambda\geq C$ and $s\geq C(T^5+T^{10})$.

 \end{lem}
 This lemma, for the weight defined in (\ref{weight}) with $m=4$, is the main result in \cite{Fer}. Again, in order to prove it one can follow the steps of the proof in \cite{Fer}, keeping in mind estimates (\ref{alpha}). \\
To finish, we give a regularity result which will be very useful for our purpose:
\begin{lem}
 Let $a\in \mathbf{R}$ and $B\in \mathbf{R}^N$ be constant and let us assume that $f\in L^2(0,T;V)$. Then, there exists a unique solution $u\in L^2(0,T;H^3(\Omega)^N\cap V)\cap H^1(0,T;V)$, together with some $p$, to the Stokes system
\begin{equation}\label{eq5}
   \begin{cases}
       u_t - \Delta u +au+B\cdot \nabla u +\nabla p=f ,\nabla \cdot u=0 &\quad \text{in} \  Q,\\
       u=0 &\quad \text{on}\  \Sigma, \\
       u_{|t=0}=0 &\quad \text{in} \ \Omega,
   \end{cases}
\end{equation}
 and there exists a constant $C>0$ such that 
\begin{equation}\label{eq6}
\|u\|_{L^2(0,T;H^3(\Omega)^N)}+\|u\|_{H^1(0,T;H^1(\Omega)^N)}\leq C\|f\|_{L^2(0,T;H^1(\Omega)^N)}.
\end{equation}
\end{lem}
This result can be found in \cite{lady}. A proof is also given in \cite{Gue}.

\section{Carleman Estimate}
In this section, we will prove a Carleman estimate which leads to an observability inequality, which in turn implies the null controllability of a linear system, similar to the linearized system associated to (\ref{2}). 
This inequality will be the main tool in the proof of Theorem $1.1$.\\   
Here, we consider the following coupled Stokes system :
\begin{equation}\label{eq7}
 \begin{cases}
       -\varphi_t - \Delta \varphi +\nabla \pi=\psi\mathbf{ 1}_{\mathcal O}+g_0,\nabla \cdot \varphi=0  &\quad \text{in} \  Q,\\
       \psi_t - \Delta \psi +\nabla \kappa=g_1,\nabla \cdot \psi=0  &\quad \text{in} \  Q,\\
       \varphi=\psi=0 &\quad \text{on}\  \Sigma, \\
       \varphi_{|t=T}=\varphi_0 ,\psi_{|t=0}=\psi_0  &\quad \text{in} \ \Omega,
   \end{cases}
\end{equation}
where $g_0,\ g_1\in L^2(Q)^N$ and $ \psi_0,\ \varphi_0 \in H$.\\
System (\ref{eq7}) is the non-homogeneous formal adjoint of the linearized of (\ref{2}) around $(0,0)$. We will be led to prove, for an open set $\omega_0\subset \mathcal O\cap\omega$, the following kind of observability inequality for (\ref{eq7}):
\begin{equation}\label{obs}
 \iint\limits_{Q}e^{-C_1/t^{m}}(|\varphi|^2+|\psi|^2)\textrm{d}x\textrm{d}t\leq C\left(\quad\iint\limits_{\omega_0\times(0,T)}|\varphi|^2\textrm{d}x\textrm{d}t+\iint\limits_{Q}e^{-C_2/t^{m}}(|g_0|^2+|g_1|^2)\textrm{d}x\textrm{d}t\right),
\end{equation}
for some $m>0$ and  certain positive constants $C$, $C_1$, $C_2$ depending on $\Omega$, $\omega_0$ and $T$ but independent of $\psi_0$ and $\varphi_0$. % and some positive numbers $m$.
 To prove such an inequality, usually, we use a combination of observability inequalities for both $\varphi$ and $\psi$ and try to eliminate the local term in $\psi$.
Even in the simpler situation of the Stokes system ($g_0\equiv g_1\equiv 0$), due to the pressure term, one cannot expect to achieve such an objective this way (see \cite{Gue}, for an explanation of this fact).\\
We will prove the following result:
\begin{theoreme} 
There exists a constant $C>0$ which depends on $\Omega$, $\omega$, $\mathcal O$ and $T$ such that
\begin{displaymath}
 \widetilde I(s,\lambda;\nabla\times\psi)+I_1(s,\lambda;\varphi)\leq C\left(s^{15}\lambda^{16}\iint\limits_{\omega\times (0,T)}e^{-4s\alpha^{\star}+s\alpha}\xi^{15}|\varphi|^2\textrm{d}x\textrm{d}t\right.
\end{displaymath}
\begin{equation}\label{eq8}
\left.+s^{5}\lambda^{6}\iint\limits_{Q}e^{-2s\alpha^{\star}-2s\alpha}\xi^{5}|g_0|^2\textrm{d}x\textrm{d}t +\iint\limits_{Q}e^{-2s\alpha^{\star}}|g_1|^2\textrm{d}x\textrm{d}t\right),
\end{equation}
for any $\lambda\geq C$, $s\geq C(T^{5}+T^{10})$, any $\varphi_0, \psi_0\in H$ and any $g_0,\ g_1\in L^2(Q)^N$, where $(\varphi, \psi)$ is the corresponding solution to (\ref{eq7}). Recall that $ \widetilde I(s,\lambda;\cdot)$ and $I_1(s,\lambda;\cdot)$ were introduced in (\ref{itilde}) and (\ref{ione}) respectively .
 \end{theoreme}
%\begin{rem}
% Combining this Carleman inequality (\ref{eq8}), with energy type estimate for both $e^{-1/t^m}\varphi$ and $e^{-1/t^m}\psi$, one can readily deduce an observability inequality like (\ref{obs}) in the case $g_0=g_1=0$ (see \cite{Gue}).
%For the proof of this theorem we will take $\varphi_{|t=T}=\varphi_0\in H$ not identically zero.
%\end{rem}

The proof of this theorem is divided in two steps. In the first step we derive a Carleman estimate for $(\nabla\times\psi)$ with a local term in $\omega_0$ using the fact that, applying the operator $(\nabla\times \cdot)$ to the second equation of system (\ref{eq7}), the resultant system can be viewed as a system of $2N-3$ heat equations. In the second one, assuming that $\psi$ is given, we apply the Carleman
estimate for Stokes systems given in Lemma 2.2. Finally we combine these two estimates and eliminate the local term in $(\nabla\times \psi)$ using the fact that $\omega_0\subset \mathcal O\cap\omega$. 
Each step will be proved in a separate paragraph.
\subsection{Carleman estimate for $\psi$}
Observe that the equation of $\psi$ is independent of $\varphi$:
\begin{equation}\label{eq9}
 \begin{cases}
       \psi_t - \Delta \psi +\nabla \kappa=g_1,\nabla \cdot \psi=0  &\quad \text{in} \  Q,\\
       \psi=0 &\quad \text{on}\  \Sigma, \\
       \psi_{|t=0}=\psi_0  &\quad \text{in} \ \Omega.
   \end{cases}
\end{equation}
A Carleman inequality for $(\nabla\times\psi)$ has been established in \cite{Gue} but for $g_1\equiv 0$. The same analysis no longer holds here since $(\nabla\times  g_1)\notin L^2(Q)^{2N-3}$. In order to get around this difficulty, we split $\psi$ (up to a weight function) into two solutions of Stokes systems. Then, we apply to the more regular one the same analysis than in \cite{Gue} and classical regularity estimates for the Stokes system to the other one.\\ 
For system (\ref{eq9}), we can prove the following result:
\begin{prop}
There exists a positive constant $C$ depending on $\Omega$ and $\omega_0$ such that 
\begin{equation}\label{eq10}
\widetilde I(s,\lambda;\nabla \times \psi)\leq C \left(s^3\lambda^4\iint\limits_{\omega_0\times (0,T)}e^{-2s\alpha}e^{-2s\alpha^{\star}}\xi^3|\nabla \times \psi|^2\textrm{d}x\textrm{d}t
+\iint\limits_{Q}e^{-2s\alpha^{\star}}|g_1|^2\textrm{d}x\textrm{d}t\right),
\end{equation}
for any $\lambda\geq C$ and  $s\geq C(T^{5}+T^{10})$. Recall that $\widetilde I(s,\lambda;\cdot)$ was defined in (\ref{itilde}).
\end{prop}

$\textbf{Proof\ of\ Proposition\ 3.1}$. Since $\psi_0\in H$ and $g_1\in L^2(Q)^N$, there exists a unique solution $(\psi,\kappa) \in L^2(0,T;V)\times D'(Q)$  of system (\ref{eq9}). Now, let $\rho(t):=e^{-s\alpha^{\star}(t)}\in C^1([0,T])$.  
Then, since $\rho$ verifies $\rho(0)=0$, $(\psi^{\star},\kappa^{\star}):=(\rho\psi,\rho\kappa)$ solves the system
\begin{equation}\label{eq11}
 \begin{cases}
       \psi^{\star}_t - \Delta \psi^{\star} +\nabla \kappa^{\star}=\rho g_1+\rho_t \psi,\nabla \cdot \psi^{\star}=0  &\quad \text{in} \  Q,\\
       \psi^{\star}=0 &\quad \text{on}\  \Sigma, \\
       \psi^{\star}_{|t=0}=0  &\quad \text{in} \ \Omega.
   \end{cases}
\end{equation}
 We decompose $(\psi^{\star},\kappa^{\star})$ as follows : $(\psi^{\star},\kappa^{\star})=(\widehat{\psi},\widehat{\kappa})+(\widetilde \psi,\widetilde{\kappa})$, where $(\widehat \psi,\widehat{\kappa})$ and $(\widetilde \psi,\widetilde{\kappa})$ solve respectively
\begin{equation}\label{eq12}
 \begin{cases}
       \widetilde{\psi}_t - \Delta  \widetilde{\psi} +\nabla \widetilde{\kappa}=\rho g_1,\nabla \cdot \widetilde{\psi} =0  &\quad \text{in} \  Q,\\
       \widetilde{\psi}=0 &\quad \text{on}\  \Sigma, \\
       \widetilde{\psi}_{|t=0}=0  &\quad \text{in} \ \Omega,
   \end{cases}
\end{equation}
and
\begin{equation}\label{eq13}
 \begin{cases}
       \widehat{\psi}_t - \Delta \widehat{\psi} +\nabla \widehat{\kappa}=\rho_t \psi,\nabla \cdot \widehat{\psi}=0  &\quad \text{in} \  Q,\\
       \widehat{\psi}=0 &\quad \text{on}\  \Sigma, \\
       \widehat{\psi}_{|t=0}=0  &\quad \text{in} \ \Omega.
   \end{cases}
\end{equation}
We apply the operator $(\nabla\times\cdot)$ to the Stokes system satisfied by $\widehat \psi$, 
 $$  (\nabla \times \widehat{\psi})_t - \Delta (\nabla \times\widehat{\psi}) =  \nabla \times(\rho_t \psi) \quad \text{in} \  Q.$$
Observe that we do not have any boundary conditions for $(\nabla \times \widehat{\psi})$. Nevertheless, we can apply Lemma $2.1$:
\begin{displaymath}
I_0(s,\lambda ;\nabla \times \widehat{\psi})\leq C\left(s^3 \lambda^4 \iint\limits_{\omega_0\times (0,T)}e^{-2s\alpha}\xi^3|\nabla \times \widehat{\psi}|^2\textrm{d}x\textrm{d}t+ 
\iint\limits_{Q}e^{-2s\alpha}(e^{-s\alpha^{\star}})_t^2|\nabla \times \psi|^2\textrm{d}x\textrm{d}t\right.
\end{displaymath}
\begin{equation}\label{eq14}
\left.+s\lambda \iint\limits_{\Sigma}e^{-2s\alpha^{\star}}\xi^{\star}\left|{\partial(\nabla \times \widehat{\psi})\over \partial n}\right|^2\textrm{d}\sigma \textrm{d}t\right),
\end{equation}
for any $\lambda\geq C$ and $s\geq C(T^{9}+T^{10}).$\\
We recall that $|\alpha^{\star}_t|\leq CT(\xi^{\star})^{6/5}$ (see (\ref{alpha})), so
\begin{equation}\label{abs}
 \iint\limits_{Q}e^{-2s\alpha}((e^{-s\alpha^{\star}})_t)^2|\nabla \times \psi|^2\textrm{d}x\textrm{d}t\leq Cs^2T^2\iint\limits_{Q}e^{-2s\alpha}e^{-2s\alpha^{\star}}(\xi^{\star})^{12/5}|\nabla \times \psi|^2\textrm{d}x\textrm{d}t,
\end{equation}
which will be absorded later on.\\
%by the left-hand side of (\ref{eq10}) 
Now, using that $\widehat{\psi}=\psi^{\star}-\widetilde \psi$ and taking into account that $(a-b)^2\geq {a^2\over2}-b^2$, we obtain:
\begin{equation}\label{car}I_0(s,\lambda ;\nabla \times \widehat{\psi})\geq {1\over 2}s^3 \lambda^4\iint\limits_{Q}e^{-2s\alpha}\xi^3|\nabla\times \psi^{\star}|^2\textrm{d}x\textrm{d}t+
{1\over 2}s\lambda^2\iint\limits_{Q}e^{-2s\alpha}\xi|\nabla (\nabla\times \psi^{\star})|^2\textrm{d}x\textrm{d}t
\end{equation}
$$-s^3 \lambda^4\iint\limits_{Q}e^{-2s\alpha}\xi^3|\nabla\times \widetilde{\psi}|^2\textrm{d}x\textrm{d}t-s\lambda^2\iint\limits_{Q}e^{-2s\alpha}\xi|\nabla (\nabla\times \widetilde{\psi})|^2\textrm{d}x\textrm{d}t.$$
Observe that the first term in the right-hand side of (\ref{car}) absorbs (\ref{abs}) as long as $\lambda\geq 1$ and $s\geq CT^{8}$.\\
We turn to the equation satisfied by $\widetilde{\psi}$. Using regularity results for system (\ref{eq12}), (see \cite{Temam}, Proposition 2.2) we deduce that
\begin{equation}\label{equ1}
 s^3 \lambda^4\iint\limits_{Q}e^{-2s\alpha}\xi^3|\nabla\times \widetilde{\psi}|^2\textrm{d}x\textrm{d}t\leq C\iint\limits_{Q}|\nabla\times\widetilde{\psi}|^2\textrm{d}x\textrm{d}t\leq C\|\widetilde{\psi}\|^2_{L^2(0,T;H^1(\Omega)^N)}\leq 
C\iint\limits_{Q}e^{-2s\alpha^{\star}}|g_1|^2\textrm{d}x\textrm{d}t,
\end{equation}
and
\begin{equation}\label{equ2}
s\lambda^2\iint\limits_{Q}e^{-2s\alpha}\xi|\nabla (\nabla\times \widetilde{\psi})|^2\textrm{d}x\textrm{d}t\leq C \iint\limits_{Q}|\nabla (\nabla\times \widetilde{\psi})|^2\textrm{d}x\textrm{d}t\leq C\|\widetilde{\psi}\|^2_{L^2(0,T;H^2(\Omega)^N)}\leq 
C\iint\limits_{Q}e^{-2s\alpha^{\star}}|g_1|^2\textrm{d}x\textrm{d}t,
\end{equation}
for $\lambda\geq C$ and $s\geq CT^{10}$, with possibly differents constants $C$. Indeed, the above constants do not depend on $T$ for $\lambda\geq C$ and $s\geq CT^{10}$.\\
The next step is to estimate the local term which appears in the right-hand side of (\ref{eq14}). Again, we put $\widehat{\psi}$ in terms of $\psi^{\star}$ and $\widetilde \psi$ :

$$s^3 \lambda^4\iint\limits_{\omega_0\times (0,T)} e^{-2s\alpha}\xi^3|\nabla \times ({\psi}^{\star}-\widetilde{\psi})|^2\textrm{d}x\textrm{d}t\leq 2s^3 \lambda^4\iint\limits_{\omega_0\times (0,T)} e^{-2s\alpha}\xi^3|\nabla \times {\psi}^{\star}|^2\textrm{d}x\textrm{d}t $$
$$+2s^3 \lambda^4\iint\limits_{\omega_0\times (0,T)} e^{-2s\alpha}\xi^3|\nabla \times \widetilde{\psi}|^2\textrm{d}x\textrm{d}t.$$
Like previously
\begin{equation}\label{equ3}
s^3 \lambda^4 \iint\limits_{\omega_0\times (0,T)}e^{-2s\alpha}\xi^3|\nabla \times \widetilde{\psi}|^2\textrm{d}x\textrm{d}t\leq s^3 \lambda^4\iint\limits_{Q} e^{-2s\alpha}\xi^3|\nabla \times \widetilde{\psi}|^2\textrm{d}x\textrm{d}t\leq
C \iint\limits_{Q}e^{-2s\alpha^{\star}}|g_1|^2\textrm{d}x\textrm{d}t.
\end{equation}
At this point combining (\ref{eq14})-(\ref{equ3}), we obtain
\begin{displaymath}
%s^3 \lambda^4\iint\limits_{Q}e^{-2s\alpha}\xi^3|\nabla\times \psi^{\star}|^2\textrm{d}x\textrm{d}t+
%s\lambda^2\iint\limits_{Q}e^{-2s\alpha}\xi|\nabla (\nabla\times \psi^{\star})|^2\textrm{d}x\textrm{d}t
\widetilde I(s,\lambda;\nabla \times \psi)\leq C\left(s^3 \lambda^4\iint\limits_{\omega_0\times (0,T)} e^{-2s\alpha}\xi^3|\nabla \times {\psi}^{\star}|^2\textrm{d}x\textrm{d}t+\iint\limits_{Q}e^{-2s\alpha^{\star}}|g_1|^2\textrm{d}x\textrm{d}t\right.
\end{displaymath}
\begin{equation}\label{e1}
\left.+
s\lambda \iint\limits_{\Sigma}e^{-2s\alpha^{\star}}\xi^{\star}\left|{\partial(\nabla \times \widehat{\psi})\over \partial n}\right|^2\textrm{d}\sigma \textrm{d}t\right),
\end{equation}
for any $\lambda\geq C$ and $s\geq C(T^{5}+T^{10})$.\\
The last step will be to eliminate the boundary term in the right-hand side of (\ref{e1}).
To this end, we introduce a function $\theta \in C^2(\overline \Omega)$ such that
\begin{equation}\label{theta}\displaystyle
 {\partial \theta\over \partial n}=1\ \text{and}\ \theta=\text{constant}\ \text{on}\ \partial\Omega.
\end{equation}
Integration by parts leads to
$$s\lambda \int\limits_0^Te^{-2s\alpha^{\star}}\xi^{\star}\left(\,\int\limits\limits_{\partial\Omega}\left|{\partial(\nabla \times \widehat{\psi})\over \partial n}\right|^2\textrm{d}\sigma\right) \textrm{d}t=s\lambda \int\limits_0^Te^{-2s\alpha^{\star}}\xi^{\star}\left(\int\limits\limits_{\Omega}\Delta(\nabla\times\widehat \psi)\nabla(\nabla\times\widehat \psi)\cdot\nabla\theta \textrm{d}x\right)\textrm{d}t$$
$$+s\lambda \int\limits_0^Te^{-2s\alpha^{\star}}\xi^{\star}\left(\,\int\limits\limits_{\Omega}(\nabla\nabla\theta\nabla(\nabla\times\widehat \psi))\nabla(\nabla\times\widehat \psi)\textrm{d}x \,\right)\textrm{d}t+{s\lambda\over 2} \int\limits_0^Te^{-2s\alpha^{\star}}\xi^{\star}\left(\int\limits\limits_{\Omega}\nabla|\nabla(\nabla\times\widehat \psi)|^2\cdot\nabla\theta \textrm{d}x\right)\textrm{d}t.$$
Thus, using Cauchy-Schwarz's inequality, the above integral can be estimate as follows
 \begin{equation}\label{esti}
s\lambda \int\limits_0^Te^{-2s\alpha^{\star}}\xi^{\star}\left(\,\int\limits_{\partial \Omega}\left|{\partial(\nabla \times \widehat{\psi})\over \partial n}\right|^2d\sigma\right) \textrm{d}t\leq Cs\lambda\int\limits_0^Te^{-2s\alpha^{\star}}\xi^{\star}\|\widehat \psi\|_{H^3(\Omega)^N}\|\widehat \psi\|_{H^2(\Omega)^N}\textrm{d}t.
\end{equation}
Thanks to the interpolation inequality $\|\widehat \psi\|_{H^2(\Omega)^N} \leq \|\widehat \psi\|^{1/2}_{H^1(\Omega)^N} \|\widehat \psi\|^{1/2}_{H^3(\Omega)^N}$, we obtain
 \begin{equation}
 s\lambda\int\limits_0^Te^{-2s\alpha^{\star}}\xi^{\star}\|\widehat \psi\|_{H^3(\Omega)^N}\|\widehat \psi\|_{H^2(\Omega)^N}\textrm{d}t\leq s\lambda\int\limits_0^Te^{-2s\alpha^{\star}}\xi^{\star}\|\widehat \psi\|^{3/2}_{H^3(\Omega)^N}\|\widehat \psi\|^{1/2}_{H^1(\Omega)^N}\textrm{d}t.
\end{equation}
Finally, using Young's inequality ($ab\leq \frac{a^p}{ p} +\frac{b^{p'}}{ p'}$ with ${1\over p}+{1\over{ p'}}=1$) for $p=4$, the task reduces to estimate
\begin{equation}\label{equ4}
 s^{5/2}\lambda\int\limits_0^Te^{-2s\alpha^{\star}}(\xi^{\star})^{5/2}\|\widehat \psi\|^{2}_{H^1(\Omega)^N}\textrm{d}t+s^{1/2}\lambda\int\limits_0^Te^{-2s\alpha^{\star}}(\xi^{\star})^{1/2}\|\widehat \psi\|^{2}_{H^3(\Omega)^N}\textrm{d}t.
\end{equation}
For the first term, thanks to the fact that $\nabla\cdot\widehat\psi(t)=0$ in $\Omega$ and $\widehat\psi=\psi^*-\widetilde\psi$, we have  
$$\|\widehat \psi(t)\|_{H^1(\Omega)^N}\leq C\|\nabla\times\widehat \psi(t)\|_{L^2(\Omega)^{2N-3}}\leq C\left( \|\nabla\times\widetilde \psi(t)\|_{L^2(\Omega)^{2N-3}}+\| \nabla\times\psi^{\star}(t)\|_{L^2(\Omega)^{2N-3}} \right).$$
The first term in the right-hand side is estimated like in (\ref{equ1}) and the second one can be absorbed by the first term in the left-hand side of (\ref{e1}), for $\lambda\geq C$ and $s\geq CT^{10}$.\\
Let us estimate now the second term in (\ref{equ4}). To this end, we introduce $(\psi^{\circ},\kappa^{\circ}):=(\eta(t) \widehat\psi,\eta(t)\widehat\kappa),$ where
$$\eta(t)=s^{1/4}\lambda^{1/2}e^{-s\alpha^{\star}}(\xi^{\star})^{1/4}\quad \text{in} \quad (0,T).$$
Then, $(\psi^{\circ},\kappa^{\circ})$ fulfills
\begin{equation}\label{eq15}
 \begin{cases}
       \psi^{\circ}_t - \Delta \psi^{\circ} +\nabla {\kappa^{\circ}}=\eta\rho_t \psi+\eta_t\widehat{\psi},\nabla \cdot \psi^{\circ}=0  &\quad \text{in} \  Q,\\
       \psi^{\circ}=0 &\quad \text{on}\  \Sigma, \\
       \psi^{\circ}_{|t=0}=0  &\quad \text{in} \ \Omega.
   \end{cases}
\end{equation}
Let us prove that the right-hand side of this system belongs to $L^2(0,T;V)$. Then, we will be able to apply Lemma 2.3.  
For the first term in the right-hand side of (\ref{eq15}), we use again that $\psi$ is a divergence-free function and we get
$$\|\eta\rho_t\psi\|_{L^2(0,T;H^1(\Omega)^N)}= \left\|{\eta\rho_t\psi^{\star}\over \rho}\right\|_{L^2(0,T;H^1(\Omega)^N)}\leq C\left\|{\eta\rho_t(\nabla\times\psi^{\star})\over \rho}\right\|_{L^2(Q)^{2N-3}}. $$
Taking into account that
$$\left|{\eta\rho_t\over \rho}\right|\leq  CT s^{5/4}\lambda^{1/2}e^{-s\alpha^{\star}}(\xi^{\star})^{6/5+1/4},$$
for any $s\geq CT^{10}$, we deduce that
\begin{equation}\label{equ5}
\|\eta\rho_t\psi\|_{L^2(0,T;H^1(\Omega)^N))}\leq CTs^{5/4}\lambda^{1/2}\|e^{-s\alpha^{\star}}(\xi^{\star})^{29/20}(\nabla \times \psi^{\star})\|_{L^2(Q)^{2N-3}}.
\end{equation}
Thus, the square of this last quantity is small with respect to the first term in the left-hand side of (\ref{e1}) by taking $\lambda\geq 1$ and $s\geq CT^{6}$.\\
We turn to the second term in the right-hand side of (\ref{eq15}). Similarly as before, we have  
$$\|\eta_t\widehat\psi\|_{L^2(0,T;H^1(\Omega)^N)}\leq C\|\eta_t(\nabla \times \widehat\psi) \|_{L^2(Q)^{2N-3}}.$$
Using again that $\widehat\psi=\psi^{\star}-\widetilde \psi$, we obtain
\begin{equation}\label{equ8}
\|\eta_t\widehat\psi\|_{L^2(0,T;H^1(\Omega)^N)}\leq\|{\eta_t(\nabla\times\psi^{\star}})\|_{L^2(Q)^{2N-3}}+\|{\eta_t(\nabla\times\widetilde{\psi}})\|_{L^2(Q)^{2N-3}},
 \end{equation}
with
$$|\eta_t|\leq CTs^{5/4}\lambda^{1/2}e^{-s\alpha^{\star}}(\xi^{\star})^{29/20},$$
for any  $s\geq CT^{10}$.\\
Therefore, the first term is estimated like (\ref{equ5}) and the second one can be estimated as in (\ref{equ1}).
%where we have used
%$
%$$\left(\left\|{\eta\rho_t\widehat\psi\over \rho}\right\|_{L^2(0,T;V)}+\left\|{\eta\rho_t\widetilde{\psi}\over \rho}\right\|_{L^2(0,T;V)}\right),$$
%$$\left\|{\eta\rho_t\widehat\psi\over \rho}\right\|_{L^2(0,T;V)}\leq\left\|{\eta\rho_t(\nabla \times\widehat\psi)\over \rho}\right\|_{L^2(Q)^N}\leq CTs^{5/4}\lambda^{1/2}\|e^{-s\alpha^{\star}}(\xi^{\star})^{5/4+1/m}(\nabla \times \widehat\psi)\|_{L^2(Q)^N}
%Thus, the square of this last quantity is bounded by the left-hand side of (\ref{eq10}), since $m>4$, by definition of $\widetilde I(s,\lambda;.)$(see (\ref{itilde})). On the other hand, with an appropriate choice of $s$ and $\lambda$ we see that
%$$\left\|{\eta\rho_t\widetilde{\psi}\over \rho}\right\|_{L^2(0,T;V)}\leq\left\|{\eta\rho_t(\nabla \times \widetilde{\psi})\over \rho}\right\|_{L^2(Q)^N}\leq C \iint\limits_{Q}e^{-2s\alpha^{\star}}|g_1|^2\textrm{d}x\textrm{d}t
%.$$
%which is also bounded by $I_0(s,\lambda;\nabla \times \widehat\psi).$
Then it follows from Lemma $2.3$, that the solution of (\ref{eq15}) satisfies $\psi^{\circ}\in L^2(0,T;H^3(\Omega)^N\cap V)$ and for all $\varepsilon>0$ there exists $C_{\varepsilon}>0$ such that 
$$\| \psi^{\circ}\|^2_{L^2(0,T;H^3(\Omega)^N)}=s^{1/2}\lambda\int\limits_0^Te^{-2s\alpha^{\star}}(\xi^{\star})^{1/2}\| \widehat\psi\|^2_{H^3(\Omega)^N}dt\leq \varepsilon  \widetilde I(s,\lambda;\nabla \times \psi)+C_{\varepsilon}\iint\limits_{Q}e^{-2s\alpha^{\star}}|g_1|^2\textrm{d}x\textrm{d}t.$$
%We already have (see $\eta_t\widehat\psi$ up there)
%$$\| \psi^{\circ}\|^2_{L^2(0,T;H^1(\Omega))}=s^{5/2}\lambda\int\limits_0^Te^{-2s\alpha^{\star}}(\xi^{\star})^{5/2}\| \widehat\psi\|^2_{H^1(\Omega)}dt\leq K I_0(s,\lambda;\nabla \times \widehat\psi).$$
%Interpolating between $L^2(0,T;H^1(\Omega))$ and $L^2(0,T;H^3(\Omega))$, we obtain
%$$\| \psi^{\circ}\|^2_{L^2(0,T;H^2(\Omega))}=s^{3/2}\lambda\int\limits_0^Te^{-2s\alpha^{\star}}(\xi^{\star})^{3/2}\| \widehat\psi\|^2_{H^1(\Omega)}dt\leq \varepsilon^{1/2} I_0(s,\lambda;\nabla \times \widehat\psi).$$
%Finally, by interpolation between $L^2(0,T;H^2(\Omega))$ and $L^2(0,T;H^3(\Omega))$, and thanks to the continuity of the trace operator, we get :
%$$s\lambda\iint\limits_{\Sigma}e^{-2s\alpha^{\star}}\xi^{\star}\left|{\partial(\nabla \times \widehat\psi)\over\partial n}\right|^2d\sigma dt\leq K's\lambda\int\limits_0^Te^{-2s\alpha^{\star}}(\xi^{\star})\| \widehat\psi\|^2_{H^{5/2}(\Omega)}dt\leq \varepsilon^{3/4} I_0(s,\lambda;\nabla \times \widehat\psi).$$
%This concludes the proof of proposition $3.1$.
This, combined with (\ref{e1}) and (\ref{esti})-(\ref{equ4}), concludes the proof of Proposition 3.1.
\subsection{Carleman estimate for $\varphi$ and conclusion} 
Here we prove Theorem $3.1$, combining the results of last section and Lemma $2.2$. 
Assuming  that $\psi$ is given, we turn to the solution of
\begin{equation}\label{eq16}
 \begin{cases}
       -\varphi_t - \Delta \varphi +\nabla \pi=\psi\mathbf{ 1}_{\mathcal O}+g_0,\nabla \cdot \varphi=0  &\quad \text{in} \  Q,\\
       \varphi=0 &\quad \text{on}\  \Sigma, \\
       \varphi_{|t=T}=\varphi_0  &\quad \text{in} \ \Omega.
   \end{cases}
\end{equation}
We choose $\pi$ such that $\int\limits_{\omega_0}\pi(t,x)dx=0$ and we apply the Carleman estimate given in Lemma $2.2$, for the weight function $5\alpha\over 2$ (instead of $\alpha$). We obtain
\begin{displaymath}\
I_1(s, \lambda;\varphi)\leq C\left(s^{16} \lambda^{40}\iint\limits_{\omega_0\times (0,T)}e^{-20s\widehat\alpha+15s\alpha^{\star}}(\widehat\xi)^{16}|\varphi|^2\textrm{d}x\textrm{d}t
+s^{15/2} \lambda^{20}\iint\limits_{{\mathcal{O}}\times (0,T)}e^{-10s\widehat\alpha+5s\alpha^{\star}}(\widehat\xi)^{15/2}|\psi|^2\textrm{d}x\textrm{d}t\right.
\end{displaymath}
\begin{equation}\label{eq17}
\left.+s^{15/2} \lambda^{20}\iint\limits_{Q}e^{-10s\widehat\alpha+5s\alpha^{\star}}(\widehat\xi)^{15/2}|g_0|^2\textrm{d}x\textrm{d}t\right),
 \end{equation}
for any $\lambda\geq C$ and  $s\geq C(T^{5}+T^{10})$, where $I_1(s,\lambda;\cdot)$ is given by (\ref{ione}).

Then, the second integral in the right-hand side of (\ref{eq17}) is bounded by $\widetilde I(s,\lambda;\nabla \times \psi)$ for a suitable choice of $\lambda\geq C$ and $s\geq CT^{10}$.\\
Indeed,
\begin{displaymath}\
s^{15/2} \lambda^{20}\iint\limits_{{\mathcal{O}}\times (0,T)}e^{-10s\widehat\alpha+5s\alpha^{\star}}(\widehat\xi)^{15/2}|\psi|^2\textrm{d}x\textrm{d}t\leq Cs^{15/2} \lambda^{20}\iint\limits_{Q}e^{-10s\widehat\alpha+5s\alpha^{\star}}(\widehat\xi)^{15/2}|\nabla \times\psi|^2\textrm{d}x\textrm{d}t\leq \epsilon\widetilde I(s,\lambda;\nabla \times \psi),
\end{displaymath}
where we have used the fact that  $\|\psi\|_{L^2(\Omega)^N}\leq C\|\nabla\times \psi\|_{L^2(\Omega)^{2N-3}}$ and also that for all $\epsilon>0$ and $M\in\mathbb{R}$, there exists $C_{\epsilon,M}>0$ such that
\begin{displaymath}
 e^{s\alpha^{\star}}\leq C_{\epsilon,M}s^M\lambda^M(\widehat\xi)^M e^{s(1+\epsilon)\widehat \alpha}, 
\end{displaymath}
 for any $\lambda\geq C$ and any $s\geq CT^{10}$.\\
Now, combining the obtained inequality with (\ref{eq10}) we get:
\begin{displaymath}
 \widetilde I(s,\lambda;\nabla \times \psi)+I_1(s, \lambda;\varphi)\leq C\left(s^{16} \lambda^{40}\iint\limits_{\omega_0\times (0,T)}e^{-20s\widehat\alpha+15s\alpha^{\star}}(\widehat\xi)^{16}|\varphi|^2\textrm{d}x\textrm{d}t\right.
 \end{displaymath}
 \begin{equation}\label{eq18}
\left.+s^3\lambda^4\iint\limits_{\omega_0\times (0,T)}e^{-2s\alpha}e^{-2s\alpha^{\star}}\xi^3|\nabla \times \psi|^2\textrm{d}x\textrm{d}t
+\iint\limits_{Q}e^{-2s\alpha^{\star}}|g_1|^2\textrm{d}x\textrm{d}t+
s^{15/2} \lambda^{20}\iint\limits_{Q}e^{-10s\widehat\alpha+5s\alpha^{\star}}(\widehat\xi)^{15/2}|g_0|^2\textrm{d}x\textrm{d}t\right),
\end{equation}
for any $\lambda\geq C$ and $s\geq C(T^{5}+T^{10})$.\\
It remains to estimate the local term in $(\nabla\times \psi)$, in terms of $\varphi$. In order to do this, we use the first equation of (\ref{eq7}), where 
the coupling term appears. Since $\omega_0 \subset \mathcal O$, we have
\begin{equation}\label{equ6}
\nabla\times \psi=-(\nabla\times \varphi)_t-\Delta(\nabla\times \varphi)-(\nabla\times g_0),  \quad \text{in} \ \omega_0\times (0,T).
 \end{equation}
Thus, replacing in the second integral in right-hand side of (\ref{eq18}), we obtain:
$$s^3\lambda^4\iint\limits_{\omega_0\times (0,T)}e^{-2s\alpha}e^{-2s\alpha^{\star}}\xi^3|\nabla \times \psi|^2\textrm{d}x\textrm{d}t=
-s^3\lambda^4\iint\limits_{\omega_0\times (0,T)}e^{-2s\alpha}e^{-2s\alpha^{\star}}\xi^3(\nabla \times \psi)(\nabla\times \varphi)_t\textrm{d}x\textrm{d}t$$
$$-s^3\lambda^4\iint\limits_{\omega_0\times (0,T)}e^{-2s\alpha}e^{-2s\alpha^{\star}}\xi^3(\nabla \times \psi)(\Delta(\nabla\times\varphi) +\nabla\times g_0)\textrm{d}x\textrm{d}t.$$
We introduce an open set $\omega_1\Subset\omega$ such that $\omega_0\Subset\omega_1$ and  a positive function $\theta\in C^2_c(\omega_1)$ such that $\theta\equiv 1\ \text{in}\ \omega_0$. Then the task turns to estimate
\begin{equation}\label{equ7}
s^3\lambda^4\iint\limits_{\omega_1\times (0,T)}\theta e^{-2s\alpha}e^{-2s\alpha^{\star}}\xi^3(\nabla \times \psi)(-(\nabla\times \varphi)_t-\Delta(\nabla\times\varphi) -\nabla\times g_0)\textrm{d}x\textrm{d}t.
\end{equation}
Performing several integration by parts, in order to get out all the derivatives of $(\nabla\times\varphi)$, we get

\begin{align}s^3\lambda^4\iint\limits_{\omega_0\times (0,T)}e^{-2s\alpha}e^{-2s\alpha^{\star}}\xi^3|\nabla \times \psi|^2&\textrm{d}x\textrm{d}t\leq
s^3\lambda^4\iint\limits_{\omega_1\times (0,T)}\theta( e^{-2s\alpha}e^{-2s\alpha^{\star}}\xi^3)_t(\nabla \times \psi)(\nabla\times \varphi)\textrm{d}x\textrm{d}t\notag\\
&+s^3\lambda^4\iint\limits_{\omega_1\times (0,T)}\theta e^{-2s\alpha}e^{-2s\alpha^{\star}}\xi^3(\nabla \times \varphi)(\nabla\times g_1)\textrm{d}x\textrm{d}t\notag\\
&-2s^3\lambda^4\iint\limits_{\omega_1\times (0,T)}\nabla (\theta e^{-2s\alpha}e^{-2s\alpha^{\star}}\xi^3)\cdot\nabla(\nabla \times \psi)(\nabla\times \varphi)\textrm{d}x\textrm{d}t\notag\\
&-s^3\lambda^4\iint\limits_{\omega_1\times (0,T)}\Delta (\theta e^{-2s\alpha}e^{-2s\alpha^{\star}}\xi^3)(\nabla \times \psi)(\nabla\times \varphi)\textrm{d}x\textrm{d}t\notag\\
&-s^3\lambda^4\iint\limits_{\omega_1\times (0,T)}\theta e^{-2s\alpha}e^{-2s\alpha^{\star}}\xi^3(\nabla \times \psi)(\nabla\times g_0)\textrm{d}x\textrm{d}t.\label{ipp}
\end{align}

Here, we have used the equation satisfied by $(\nabla\times\psi)$ and the fact that $\theta$ has compact support in $\omega_1$.
We perform another integration by parts and use Young's inequality to obtain:
$$s^3\lambda^4\iint\limits_{\omega_1\times (0,T)}\theta e^{-2s\alpha}e^{-2s\alpha^{\star}}\xi^3(\nabla \times \varphi)(\nabla\times g_1)\textrm{d}x\textrm{d}t=
-s^3\lambda^4\iint\limits_{\omega_1\times (0,T)}\nabla\times (\theta e^{-2s\alpha}e^{-2s\alpha^{\star}}\xi^3(\nabla \times \varphi) )g_1\textrm{d}x\textrm{d}t$$
$$\leq C\left( \iint\limits_{Q} e^{-2s\alpha^{\star}}|g_1|^2\textrm{d}x\textrm{d}t + s^6\lambda^8\iint\limits_{Q}e^{-4s\alpha}e^{-2s\alpha^{\star}}\xi^6(s^2\lambda^2\xi^2|\nabla\times \varphi|^2
+|\nabla(\nabla\times \varphi)|^2)\textrm{d}x\textrm{d}t\right).$$
The last term in this inequality is estimated by $\epsilon I_1(s,\lambda;\varphi)$ for $\lambda\geq C$ and $s\geq CT^{10}$. An analogous estimate holds for the term containing $(\nabla \times g_0)$: %$ e^{-4s\alpha-7s\alpha^{\star}+6s\widehat\alpha}$ is  smaller than $e^{-2s\alpha-2s\alpha^{\star}}$ for $\lambda\geq C$ and $s\geq CT^{10}$.\\
$$s^3\lambda^4\iint\limits_{\omega_1\times (0,T)}\theta e^{-2s\alpha}e^{-2s\alpha^{\star}}\xi^3(\nabla \times \psi)(\nabla\times g_0)\textrm{d}x\textrm{d}t=
-s^3\lambda^4\iint\limits_{\omega_1\times (0,T)}\nabla\times (\theta e^{-2s\alpha}e^{-2s\alpha^{\star}}\xi^3(\nabla \times \psi) )g_0\textrm{d}x\textrm{d}t$$
$$\leq Cs^5\lambda^6 \iint\limits_{Q}\xi^5 e^{-2s\alpha-2s\alpha^{\star}}|g_0|^2\textrm{d}x\textrm{d}t +\epsilon \widetilde I(s,\lambda;\nabla \times \psi).$$
On the other hand we have the following estimates for the weight functions:
$$|( e^{-2s\alpha}e^{-2s\alpha^{\star}}\xi^3)_t|\leq CTse^{-2s\alpha}e^{-2s\alpha^{\star}}(\xi)^{4+1/5}\quad \text{and} \quad
|\Delta (e^{-2s\alpha}e^{-2s\alpha^{\star}}\xi^3)|\leq Cs^2\lambda^2e^{-2s\alpha}e^{-2s\alpha^{\star}}\xi^5,$$
for any $s\geq CT^{10}$.\\
Using these estimates for the first, third and fourth terms in the right-hand side of (\ref{ipp}), we deduce that
$$s^3\lambda^4\iint\limits_{\omega_0\times (0,T)}e^{-2s\alpha}e^{-2s\alpha^{\star}}\xi^3|\nabla \times \psi|^2\textrm{d}x\textrm{d}t\leq \epsilon (\widetilde I (s,\lambda;\nabla \times \psi)+ I_1(s,\lambda;\varphi))
$$
\begin{equation}\label{etoile} +C_{\epsilon}\left(s^7\lambda^8 \iint\limits_{\omega_1\times (0,T)}e^{-2s\alpha}e^{-2s\alpha^{\star}}\xi^7|\nabla \times \varphi|^2\textrm{d}x\textrm{d}t+\iint\limits_{\omega_1\times (0,T)}e^{-2s\alpha^{\star}}|g_1|^2\textrm{d}x\textrm{d}t+
s^5\lambda^6\iint\limits_{Q}e^{-2s\alpha-s\alpha^{\star}}\xi^5|g_0|^2\textrm{d}x\textrm{d}t\right), 
\end{equation}
for $\lambda\geq C$ and $s\geq C(T^{5}+T^{10})$.\\
Furthermore, considering an open set $\omega_2\Subset \omega$ such that $\omega_1\Subset \omega_2$, one can prove that
$$s^7\lambda^8 \iint\limits_{\omega_1\times (0,T)}e^{-2s\alpha }e^{-2s\alpha^{\star}}\xi^7|\nabla \times \varphi|^2\textrm{d}x\textrm{d}t\leq 
\epsilon s^{-1}\iint\limits_{Q}e^{-5s\alpha }\xi^{-1}|\Delta \varphi|^2\textrm{d}x\textrm{d}t+s^{15}\lambda^{16}\iint\limits_{\omega_2\times (0,T)}e^{-4s\alpha^{\star}+s\alpha }\xi^{15}|\varphi|^2\textrm{d}x\textrm{d}t.$$
This, combined with (\ref{etoile}) and (\ref{eq18}), gives the desired inequality (\ref{eq8}).

\section{Null controllability of the linear system}
In this section we consider a linear coupled Stokes system with right-hand sides. More precisely, we look for a control $v\in L^2(\omega\times (0,T))^N$ such that, under suitable decreasing properties on $f_1$ and $f_2$, the solution to
\begin{equation} \label{eq19} \left\lbrace { \begin{tabular}{llll}
 $w_t- \Delta w + \nabla p=f_1+v\mathbf{ 1}_{\omega},\ \nabla \cdot w=0 $ & \textrm{in }$Q$,\\
 $-z_t- \Delta z + \nabla q=f_2+w\mathbf{ 1}_{\mathcal O},\ \nabla \cdot z=0 $ & \textrm{in }$Q$,\\
$w=z=0$ & \textrm{on }$\Sigma$,\\
$  w_{|t=0}=z_{|t=T}=0$ & \textrm{in }$\Omega$,\\ \end{tabular}
} \right. 
\end{equation}
satisfies
\begin{equation}\label{eq20}
z_{|t=0}=0 \quad \text{in} \quad \Omega.
\end{equation}
As we have already mentioned, an observability inequality for (\ref{eq7}) will imply the null controllability of (\ref{eq19}) with decreasing properties for the state(s) and the control(s) (see (\cite{Fer})).
Here, we present a null controllability result for (\ref{eq19}) where we look for a more regular solution $(w,z)$. This will be done by solving the controllability problem in  spaces depending on the previous weight functions. Furthermore, this result will be useful to deduce the local null controllability of the nonlinear problem (\ref{2}) in the last section.\\
First let us prove a modified Carleman inequality, from (\ref{eq8}), with weight functions that do not vanish at $t=T$. To be more specific, consider
\begin{equation}\label{l(t)}
l(t)=\begin{cases}
t(T-t),\ 0\leq t\leq T/2,\\
\displaystyle{T^{2}\over 4},\ T/2\leq t\leq T,
\end{cases} 
 \end{equation}
and the following associated weight functions :
\begin{equation}
 \beta(x,t)={\exp(12\lambda \|\eta^0\|_{\infty})-\exp{\lambda(k\|\eta^0\|_{\infty}+\eta^0(x))}\over l(t)^5},\quad\gamma(x,t)={\exp{\lambda(10\|\eta^0\|_{\infty}+\eta^0(x))}\over l(t)^5}
\end{equation}

\begin{equation}
\beta^{\star}(t)=\max_{x\in \overline \Omega}\beta(x,t),\quad  \widehat\beta(t)=\min_{x\in \overline \Omega}\beta(x,t),\quad\gamma^{\star}(t)=\min_{x\in \overline \Omega}\gamma(x,t), \quad \widehat\gamma(t)=\max_{x\in \overline \Omega}\gamma(x,t).
\end{equation}
With this definition we have the following

%\begin{lem}\label{Carleman}}
%Let $s$ and $\lambda$ like in Theorem $3.1$. Then, there exists a positive constant $C$ depending on $\Omega$, $\omega_0$, $T$, $s$ and $\lambda$ such that
%\begin{equation} \label{22}
%\iint\limits_{Q}e^{-4s\beta^{\star}}(\gamma^{\star})^3|\psi|^2\textrm{d}x\textrm{d}t+\iint\limits_{Q}e^{-5s\beta}\gamma^3|\varphi|^2\textrm{d}x\textrm{d}t\leq C\left(\iint\limits_{\omega\times (0,T)}e^{-20s\widehat\beta+15s\beta^{\star}}(\widehat\gamma)^{16}|\varphi|^2\textrm{d}x\textrm{d}t\right.
%\end{equation}
%\begin{displaymath}
%\left. +\iint\limits_{Q}e^{-10s\widehat\beta+5s\beta^{\star}}(\widehat\gamma)^{15/2}|g_0|^2\textrm{d}x\textrm{d}t +\iint\limits_{Q}e^{-2s\beta^{\star}}|g_1|^2\textrm{d}x\textrm{d}t\right),
%\end{displaymath}
%for any $\psi_0\in H$, where $(\varphi,\psi)$ is the solution of (\ref{eq7}) .
%\end{lem}

\begin{lem}\label{Carleman}
 Let $s$ and $\lambda$ like in Theorem $3.1$. Then, there exists a positive constant $C$ depending on $\Omega$, $\omega_0$, $T$, $s$ and $\lambda$ such that
\begin{equation} \label{eq22}
\iint\limits_{Q}e^{-4s\beta^{\star}}(\gamma^{\star})^3|\psi|^2\textrm{d}x\textrm{d}t+\iint\limits_{Q}e^{-5s\beta^{\star}}(\gamma^{\star})^3|\varphi|^2\textrm{d}x\textrm{d}t\leq C\left(\quad\iint\limits_{\omega\times (0,T)}e^{-4s\beta^{\star}+s\beta}\gamma^{15}|\varphi|^2\textrm{d}x\textrm{d}t\right.
\end{equation}
\begin{displaymath}
\left. +\iint\limits_{Q}e^{-2s\beta-2s\beta^{\star}}\gamma^{5}|g_0|^2\textrm{d}x\textrm{d}t +\iint\limits_{Q}e^{-2s\beta^{\star}}|g_1|^2\textrm{d}x\textrm{d}t\right),
\end{displaymath}
for any $\varphi_0,\ \psi_0\in H$, where $(\varphi,\psi)$ is the associated solution to (\ref{eq7}).
\end{lem}

$\textbf{Proof\ of\ Lemma\ \ref{Carleman}.}$
First by construction  $\alpha=\beta$ and $\xi=\gamma$ in $\Omega\times (0,T/2)$, so that
$$\int\limits_0^{T/2}\!\!\int\limits\limits_{\Omega}e^{-4s\alpha^{\star}}(\xi^{\star})^3|\psi|^2\textrm{d}x\textrm{d}t+\int\limits_0^{T/2}\!\!\int\limits\limits_{\Omega}e^{-5s\alpha^{\star}}(\xi^{\star})^3|\varphi|^2\textrm{d}x\textrm{d}t$$
$$=\int\limits_0^{T/2}\!\!\int\limits\limits_{\Omega}e^{-4s\beta^{\star}}(\gamma^{\star})^3|\psi|^2\textrm{d}x\textrm{d}t+\int\limits_0^{T/2}\!\!\int\limits\limits_{\Omega}e^{-5s\beta^{\star}}(\gamma^{\star})^3|\varphi|^2\textrm{d}x\textrm{d}t.$$
Therefore, it follows from (\ref{eq8}) (observe that $e^{-4s\beta^{\star}}(\gamma^{\star})^3\leq e^{-2s\beta}e^{-2s\beta^{\star}}\gamma^3$, $e^{-5s\beta^{\star}}(\gamma^{\star})^3\leq e^{-5s\beta}\gamma^3$ and $\|\psi\|_{L^2(\Omega)^N}\leq C\|\nabla\times \psi\|_{L^2(\Omega)^{2N-3}}$)
$$\int\limits_0^{T/2}\!\!\int\limits\limits_{\Omega}e^{-4s\beta^{\star}}(\gamma^{\star})^3|\psi|^2\textrm{d}x\textrm{d}t+\int\limits_0^{T/2}\!\!\int\limits\limits_{\Omega}e^{-5s\beta^{\star}}(\gamma^{\star})^3|\varphi|^2\textrm{d}x\textrm{d}t\leq C(T,s,\lambda)\left(\quad
\iint\limits_{\omega\times (0,T)}e^{-4s\alpha^{\star}+s\alpha}\xi^{15}|\varphi|^2\textrm{d}x\textrm{d}t\right.$$
$$\left.+\iint\limits_{Q}e^{-2s\alpha^{\star}}|g_1|^2\textrm{d}x\textrm{d}t+\iint\limits_{Q}e^{-2s\alpha-2s\alpha^{\star}}\xi^{5}|g_0|^2\textrm{d}x\textrm{d}t\right),$$
for any $\psi_0\in H$.\\
Thus, by definition of $\beta,\beta^{\star},\gamma$ and $\gamma^{\star}$ we have
\begin{displaymath}
 \int\limits_0^{T/2}\!\!\int\limits\limits_{\Omega}e^{-4s\beta^{\star}}(\gamma^{\star})^3|\psi|^2\textrm{d}x\textrm{d}t+\int\limits_0^{T/2}\!\!\int\limits\limits_{\Omega}e^{-5s\beta^{\star}}(\gamma^{\star})^3|\varphi|^2\textrm{d}x\textrm{d}t\leq C(T,s,\lambda)\left(\quad
 \iint\limits_{\omega\times (0,T)}e^{-4s\beta^{\star}+s\beta}\gamma^{15}|\varphi|^2\textrm{d}x\textrm{d}t\right.
\end{displaymath}
\begin{equation}\label{0t2}
\left.+\iint\limits_{Q}e^{-2s\beta^{\star}}|g_1|^2\textrm{d}x\textrm{d}t+\iint\limits_{Q}e^{-2s\beta-2s\beta^{\star}}\gamma^{5}|g_0|^2\textrm{d}x\textrm{d}t\right),
\end{equation}

We turn to the domain $\Omega\times (T/2,T)$. Here, we will use  well known a priori estimates for the Stokes system.
Indeed, let us introduce a function  $\zeta\in C^1([0,T])$ such that
$$ \zeta=0\ \text{in}\ [0,T/4],\quad \zeta=1\ \text{in}\ [T/2,T],\quad |\zeta'|\leq C/T.$$
Using classical energy estimates for both $\zeta\varphi$ and $\zeta\psi$ (see, for instance, \cite{lady}), which solve the Stokes system (\ref{eq7}), we obtain
$$\|\zeta\varphi\|^2_{L^2(T/4,T;H^1(\Omega)^N)} +\|\zeta\varphi\|^2_{L^{\infty}(T/4,T;H)}\leq C\left( \|\zeta g_0\|^2_{L^2(T/4,T;L^2(\Omega)^N)}+ \|\zeta\psi\|^2_{L^2(T/4,T;L^2(\mathcal O)^N)}\right.$$
$$\left.+ {1\over T^2}\|\varphi\|^2_{L^2(T/4,T/2;L^2(\Omega)^N)} \right)$$
and
$$\|\zeta\psi\|^2_{L^2(T/4,T;H^1(\Omega)^N)} +\|\zeta\psi\|^2_{L^{\infty}(T/4,T;H)}\leq C\left( \|\zeta g_1\|^2_{L^2(T/4,T;L^2(\Omega)^N)}+ {1\over T^2}\|\psi\|^2_{L^2(T/4,T/2;L^2(\Omega)^N)} \right).$$
Combining these last two inequalities and keeping in mind the definition of $\zeta$, we obtain
$$\|\varphi\|^2_{L^2(T/2,T;L^2(\Omega)^N)}+\|\psi\|^2_{L^2(T/2,T;L^2(\Omega)^N)}\leq \left(\|\zeta g_0\|^2_{L^2(T/4,T;L^2(\Omega)^N)}+ \|\zeta g_1\|^2_{L^2(T/4,T;L^2(\Omega)^N)}\right.$$
$$\left.+{1\over T^2}\|\psi\|^2_{L^2(T/4,T;L^2(\Omega)^N)}+ {1\over T^2}\|\varphi\|^2_{L^2(T/4,T/2;L^2(\Omega)^N)}\right)$$
Using (\ref{0t2}) to estimate the last two terms and taking into account that the weight functions $\beta$ and $\gamma$ are bounded in $[T/4,T]$, we get the following estimate
\begin{equation}\label{t2t}
 \int\limits_{T/2}^T\!\!\int\limits\limits_{\Omega}e^{-4s\beta^{\star}}(\gamma^{\star})^3|\psi|^2\textrm{d}x\textrm{d}t+\int\limits_{T/2}^T\!\!\int\limits_{\Omega}e^{-5s\beta^{\star}}(\gamma^{\star})^3|\varphi|^2\textrm{d}x\textrm{d}t\leq C(T,s,\lambda)\left(\, \int\limits_{T/4}^T\!\!\int\limits\limits_{Q}e^{-2s\beta^{\star}}|g_1|^2\textrm{d}x\textrm{d}t\right.
\end{equation}
\begin{displaymath}
\left.+\int\limits_{T/4}^T\!\!\int\limits_{Q}e^{-2s\beta-2s\beta^{\star}}\gamma^{5}|g_0|^2\textrm{d}x\textrm{d}t\right).
\end{displaymath}
This, together with (\ref{0t2}), gives us the desired inequality (\ref{eq22}).

\vskip0.5cm Now, we will use this Carleman inequality to deduce a null controllability result for system (\ref{eq19}). In the same spirit of \cite{Fer}, where the local exact controllability of the Navier-Stokes system is proved, we introduce the following weighted space:
\begin{displaymath}
 \mathcal E^{s,\lambda}=\left\{(w,z,p,q,v):e^{s\beta+s\beta^{\star}}(\gamma^{\star})^{-5/2}w\in L^2(Q)^N,\  e^{s\beta^{\star}}z\in L^2(Q)^N,\ 
e^{2s\beta^{\star}-{1\over2}s\beta}\gamma^{-15/2}v\mathbf{1}_{\omega}\in L^2(Q)^N,\right.\   
\end{displaymath}
\begin{equation}\label{eq21}   
\left.e^{{5\over2}s\beta^{\star}}(\gamma^{\star})^{-3/2}(w_t- \Delta w + \nabla p-v\mathbf{1}_{\omega})\in L^2(Q)^N,\ 
 e^{2s\beta^{\star}}(\gamma^{\star})^{-3/2}(-z_t- \Delta z + \nabla q-w\mathbf{ 1}_{\mathcal O})\in L^2(Q)^N\right\}. 
\end{equation}
Defined as we have seen, $\mathcal E^{s,\lambda}$ is a Banach space for the norm
$$ \left\|(w,z,p,q,v)\right\|_{ \mathcal E^{s,\lambda}}=\left(\left\|e^{s\beta+s\beta^{\star}}(\gamma^{\star})^{-5/2}w\right\|^2_{L^2(Q)^N}+\left\|e^{s\beta^{\star}}z\right\|^2_{L^2(Q)^N}+\left\|e^{2s\beta^{\star}-{1\over2}s\beta}\gamma^{-15/2}v\mathbf{1}_{\omega}\right\|^2_{L^2(Q)^N}\right.$$
$$\left.+\left\|e^{{5\over2}s\beta^{\star}}(\gamma^{\star})^{-3/2}(w_t- \Delta w + \nabla p-v\mathbf{1}_{\omega})\right\|^2_{L^2(Q)^N}+\left\|e^{2s\beta^{\star}}(\gamma^{\star})^{-3/2}(-z_t- \Delta z + \nabla q-w\mathbf{ 1}_{\mathcal O})\right\|^2_{L^2(Q)^N}\right)^{1/2}.$$
\begin{rem}
 If $(w,z,p,q,v)\in  \mathcal E^{s,\lambda}$, then $z_{|t=0}=0$.
\end{rem}
We will prove the following result :
\begin{prop} 
 Let $f_1$, $f_2$ satisfy $e^{{5\over2}s\beta^{\star}}(\gamma^{\star})^{-3/2}f_1\in L^2(Q)^N$ and $e^{2s\beta^{\star}}(\gamma^{\star})^{-3/2}f_2\in L^2(Q)^N$. Then, there exists $v\in L^2(\omega\times (0,T))^N$ such that, if $(w,z,p,q)$ is the solution
of (\ref{eq19}), one has  $(w,z,p,q,v)\in \mathcal E^{s,\lambda}$.
\end{prop}
$\textbf{Proof\ of\ Proposition\ 4.1.}$
Let us introduce the following constrained extremal problem:

\begin{equation}\label{eq23}
 \begin{cases}\displaystyle\inf {1\over 2}\left(\iint\limits_{Q}e^{2s\beta+2s\beta^{\star}}\gamma^{-5}|w|^2\textrm{d}x\textrm{d}t+\iint\limits_{Q}e^{2s\beta^{\star}}|z|^2\textrm{d}x\textrm{d}t+\iint\limits_{\omega\times (0,T)}e^{4s\beta^{\star}-s\beta}\gamma^{-15}|v|^2\textrm{d}x\textrm{d}t\right)\\
\text{subject to } v\in L^2(Q)^N,\  \text{supp}\ v\subset \omega\times (0,T)\ \text{and}\\
\left\lbrace {\begin{tabular}{llll}
 $w_t- \Delta w + \nabla p=f_1+v\mathbf{ 1}_{\omega},\ \nabla \cdot w=0 $ & \textrm{in }$Q$,\\
 $-z_t- \Delta z + \nabla q=f_2+w\mathbf{ 1}_{\mathcal O},\ \nabla \cdot z=0 $ & \textrm{in }$Q$,\\
$w=z=0$ & \textrm{on }$\Sigma$,\\
$  w_{|t=0}=z_{|t=T}=z_{|t=0}=0$ & \textrm{in }$\Omega$. \end{tabular}
} \right. 
\end{cases}
\end{equation}
%\begin{equation}\label{optcon}
% \begin{cases}
 %\displaystyle \inf {1\over 2}\left(\int\limits\!\!\!\!\int\limits\limits_{Q}e^{6s\widehat\beta-3s\beta^{\star}}(\widehat\gamma)^{-15/2}|w|^2\textrm{d}x\textrm{d}t+\int\limits\!\!\!\!\int\limits\limits_{Q}e^{2s\beta^{\star}}|z|^2\textrm{d}x\textrm{d}t+\int\limits\!\!\!\!\int\limits_{\omega\times (0,T)}e^{12s\widehat\beta-9s\beta^{\star}}(\widehat\gamma)^{-16}|v|^2\textrm{d}x\textrm{d}t\right)\\
 %\
%\text{subject to } v\in L^2(Q)^N,\  \text{supp}\ v\subset \omega\times (0,T)\ \text{and}\\
%\left\lbrace { \begin{array}{llll}
 %$w_t- \Delta w + \nabla p=f_1+v\mathbf{ 1}_{\omega},\ \nabla \cdot w=0 $ & \textrm{in }$Q$,\\
 %$-z_t- \Delta z + \nabla q=f_2+w\mathbf{ 1}_{\mathcal O},\ \nabla \cdot z=0 $ & \textrm{in }$Q$,\\
%$w=z=0$ & \textrm{on }$\Sigma$,\\
%$  w_{|t=0}=z_{|t=T}=z_{|t=0}=0$ & \textrm{in }$\Omega$.\\ \end{array}
%} \right. 
%\end{cases}
%\end{equation}

Assume that this problem admits a unique solution $(\widehat w,\widehat z,\widehat p,\widehat q,\widehat v)$. Then, in view of the Lagrange's principle there exists dual variables $(\overline{w},\overline{z},\overline{p},\overline{q})$ such that

\begin{equation}\label{eq24}
 \begin{cases}
\widehat w=e^{-2s\beta-2s\beta^{\star}}\gamma^{5}(-\overline{w}_t- \Delta \overline{w} + \nabla \overline{p}-\overline z\mathbf{ 1}_{\mathcal O}) &\quad \text{in} \  Q,\\
\widehat z=e^{-2s\beta^{\star}}(\overline{z}_t- \Delta \overline{z} + \nabla \overline{q}) &\quad \text{in} \  Q,\\
\widehat v=e^{-4s\beta^{\star}+s\beta}\gamma^{15}\overline{w} &\quad \text{in} \  \omega\times (0,T),\\
\widehat w=\widehat z =0 &\quad \text{on}\  \Sigma.
\end{cases}
\end{equation}

Let us set 
$$\mathcal P_0=\{(w,z,p,q)\in C^{\infty}(\overline Q)^{2N+2}\ ;\nabla \cdot w=\nabla \cdot z=0\ \text{in}\ Q , w=z=0\ \text{on}\ \Sigma\ \text{and}\ \int\limits_{\omega_0}q(x,t)dx=0 \}$$
and	
\begin{equation}\label{eq28}	
a((\overline w,\overline z,\overline p,\overline q),(w,z,p,q))=\iint\limits_{Q}e^{-2s\beta-2s\beta^{\star}}\gamma^{5}(-\overline{w}_t- \Delta \overline{w} + \nabla \overline{p}-\overline z\mathbf{ 1}_{\mathcal O}) (-w_t- \Delta w+ \nabla p-z\mathbf{ 1}_{\mathcal O}) \textrm{d}x\textrm{d}t
\end{equation}
\begin{displaymath}
 +\iint\limits_{Q}e^{-2s\beta^{\star}}(\overline{z}_t- \Delta \overline{z} + \nabla \overline{q})(z_t- \Delta z + \nabla q)\textrm{d}x\textrm{d}t+\iint\limits_{\omega\times (0,T)}e^{-4s\beta^{\star}+s\beta}\gamma^{15}\overline w w\textrm{d}x\textrm{d}t\quad \forall (w,z,p,q)\in \mathcal P_0.
\end{displaymath}
With this definition, one can see that, if the functions $\widehat w$, $\widehat z$ and $\widehat v$ solve (\ref{eq23}), we must have
\begin{equation}\label{eq25}
 a((\overline w,\overline z,\overline p,\overline q),(w,z,p,q))= l(w,z,p,q),\quad \forall (w,z,p,q)\in P_0,
\end{equation}
where 
\begin{equation}\label{eq26}
l(w,z,p,q)=\iint\limits_{Q}f_1w\textrm{d}x\textrm{d}t+\iint\limits_{Q}f_2z\textrm{d}x\textrm{d}t.
\end{equation}
The main idea is to prove that there exists exactly one $(\overline{w},\overline{z},\overline{p},\overline{q})$ satisfying (\ref{eq25}). Then we will define $(\widehat w,\widehat z,\widehat p,\widehat q,\widehat v)$ using (\ref{eq24}) and we will check that it fulfills the 
desired properties.\\
Indeed, observe that the Carleman inequality (\ref{eq22}) holds for  $(w,z,p,q)\in \mathcal P_0$,
\begin{equation}\label{eq27}
\iint\limits_{Q}e^{-4s\beta^{\star}}(\gamma^{\star})^3|w|^2\textrm{d}x\textrm{d}t+\iint\limits_{Q}e^{-5s\beta^{\star}}(\gamma^{\star})^3|z|^2\textrm{d}x\textrm{d}t\leq C a(( w, z, p, q),(w,z,p,q))\quad \forall (w,z,p,q)\in \mathcal P_0.
\end{equation}

In the linear space $\mathcal P_0$ we consider the bilinear form $a(.,.)$ given by (\ref{eq28}); from the unique continuation property for Stokes-like systems (see \cite{Fabre}) %\begin{equation}\label{eq29}
 %  \begin{cases}
  %     u_t - \Delta u +\nabla h=f ,\nabla \cdot u=0 &\quad \text{in} \  Q,\\
   %    u=0 &\quad \text{on}\  \Sigma, \\
    %   u_{|t=0}=u^0 &\quad \text{in} \ \Omega,
   %\end{cases}
%\end{equation}
we deduce that $a(.,.)$ is a scalar product in $\mathcal P_0$.
Let us now consider the space $\mathcal P$, given by the completion of $\mathcal P_0$ for the norm associated to $a(.,.)$. This is a Hilbert space and $a(.,.)$ is a continuous and coercive bilinear form on $\mathcal P$.\\
We turn to the linear operator  $l$, given by (\ref{eq26}) for all $(w,z,p,q)\in \mathcal P$, a simple computation leads to 
$$ l(w,z,p,q)\leq \|e^{{5\over2}s\beta^{\star}}(\gamma^{\star})^{-3/2}f_1\|_{L^2(Q)^N}\|e^{-{5\over2}s\beta^{\star}}(\gamma^{\star})^{3/2}w\|_{L^2(Q)^N}+\| e^{2s\beta^{\star}}(\gamma^{\star})^{-3/2}f_2\|_{L^2(Q)^N}\| e^{-2s\beta^{\star}}(\gamma^{\star})^{3/2}z\|_{L^2(Q)^N}$$
Then, using (\ref{eq27}) and the density of $\mathcal P_0$ in $\mathcal P$, we have 
$$l(w,z,p,q)\leq C (\|e^{{5\over2}s\beta^{\star}}(\gamma^{\star})^{-3/2}f_1\|_{L^2(Q)^N}+\|e^{2s\beta^{\star}}(\gamma^{\star})^{-3/2}f_2\|_{L^2(Q)^N})\|(w,z,p,q)\|_{\mathcal P} \quad \forall (w,z,p,q)\in \mathcal P.$$
Consequently $l$ is a bounded linear operator on $\mathcal P$. Then, in view of Lax-Milgram's lemma, there exists one and only one $(\overline w,\overline z,\overline p,\overline q)$ satisfying

\begin{equation}\label{eq30}
 \begin{cases}
a((\overline w,\overline z,\overline p,\overline q),(w,z,p,q))= l(w,z,p,q),\quad \forall (w,z,p,q)\in \mathcal P\\
(\overline w,\overline z,\overline p,\overline q)\in P.
\end{cases}
\end{equation}
We finally get the existence of $(\widehat w,\widehat z,\widehat p,\widehat q,\widehat v)$, just setting
$$\widehat w=e^{-2s\beta-2s\beta^{\star}}\gamma^{5}(-\overline{w}_t- \Delta \overline{w} + \nabla \overline{p}-\overline z\mathbf{ 1}_{\mathcal O}) ,\quad \widehat z=e^{-2s\beta^{\star}}(\overline{z}_t- \Delta \overline{z} + \nabla \overline{q})
\ \text{and}\quad
\widehat v=e^{-4s\beta^{\star}+s\beta}\gamma^{15}\overline{w} .$$\\
It remains to check that $(\widehat w,\widehat z,\widehat p,\widehat q,\widehat v)$ verifies
$$\iint\limits_{Q}e^{2s\beta+2s\beta^{\star}}\gamma^{-5}|\widehat w|^2\textrm{d}x\textrm{d}t+\iint\limits_{Q}e^{2s\beta^{\star}}|\widehat z|^2\textrm{d}x\textrm{d}t+\iint\limits_{\omega\times (0,T)}e^{4s\beta^{\star}-s\beta}\gamma^{-15}|\widehat v|^2\textrm{d}x\textrm{d}t< +\infty$$
and solves the Stokes system in (\ref{eq23}). 
The first point is easy to check, since $(\overline w,\overline z,\overline p,\overline q)\in \mathcal P$ and 
$$\iint\limits_{Q}e^{2s\beta+2s\beta^{\star}}\gamma^{-5}|\widehat w|^2\textrm{d}x\textrm{d}t+\iint\limits_{Q}e^{2s\beta^{\star}}| \widehat z|^2\textrm{d}x\textrm{d}t+\iint\limits_{\omega\times (0,T)}e^{4s\beta^{\star}-s\beta}\gamma^{-15}|\widehat v|^2\textrm{d}x\textrm{d}t=a((\overline w,\overline z,\overline p,\overline q),(\overline w,\overline z,\overline p,\overline q))< +\infty.$$
In order to check the second point, we introduce the (weak) solution $(\widetilde w,\widetilde z,\widetilde p,\widetilde q)$  to the Stokes system
\begin{equation}\label{eq31}
\begin{cases}
 \widetilde{w}_t- \Delta \widetilde{w} + \nabla \widetilde p=f_1+\widehat v\mathbf{ 1}_{\omega},\ \nabla \cdot \widetilde w=0  & \textrm{in }Q,\\
 -\widetilde{z}_t- \Delta \widetilde{z} + \nabla \widetilde q=f_2+\widetilde w\mathbf{ 1}_{\mathcal O},\ \nabla \cdot \widetilde z=0  & \textrm{in }Q,\\
  \widetilde w=\widetilde z=0 & \textrm{on }\Sigma,\\
 \widetilde{w}_{|t=0}=\widetilde{z}_{|t=T}=0 & \textrm{in }\Omega.
\end{cases}
\end{equation}
In particular, $(\widetilde w,\widetilde z)$ is the unique solution by transposition of (\ref{eq31}), in the following sense 
\begin{equation}\label{eq32}
 \langle(\widetilde{w},\widetilde{z}),(a,b)\rangle_{L^2(Q)^{2N}}=\langle(f_1+\widehat v\mathbf{ 1}_{\omega},f_2),(\varphi,\psi)\rangle_{L^2(Q)^{2N}},\quad \forall (a,b)\in L^2(Q)^{2N},
\end{equation}
where $(\varphi,\psi)$, together with some $(\pi,\kappa)$, solves 

\begin{equation} 
\begin{cases}
 P^{\star}(x,t;D)(\varphi,\psi)=(a,b)\quad & \textrm{in } Q,\\
\nabla\cdot\varphi=\nabla\cdot\psi=0\ \textrm{in } Q,\  \varphi=\psi=0\ & \textrm{on } \Sigma.
\end{cases}
\end{equation}

Here, we have denoted by $P^{\star}(x,t;D)$ the formal adjoint operator of $P(x,t;D)$ given by
$$P(x,t;D)(\widetilde{w},\widetilde{z})=(\widetilde{w}_t- \Delta \widetilde{w} + \nabla \widetilde p,-\widetilde{z}_t- \Delta \widetilde{z} + \nabla \widetilde q-\widetilde w\mathbf{ 1}_{\mathcal O})^t.$$
%with homogeneous Dirichlet condition, free-divergence conditions on both variables, that is to say, $(\varphi,\psi)$ solves (\ref{eq7}) with $g_0=a$ and $g_1=b$.
%Here, $\langle.,.\rangle$ is defined as follows
%$$\langle(a,b),(c,d)\rangle=\iint\limits_{Q}(a.c+b.d)\textrm{d}x\textrm{d}t\quad \forall (a,b),(c,d)\in (L^2(Q)^N)^2.$$
From (\ref{eq30}) and the definition of  $(\widehat w,\widehat z,\widehat v)$, we see that $(\widehat w,\widehat z)$ also satisfies (\ref{eq32}). Consequently, $(\widehat w,\widehat z)=(\widetilde{w},\widetilde{z})$ and $(\widehat w,\widehat z,\widehat p,\widehat q)$ is the solution to the Stokes system (\ref{eq23}).

\section{Insensitizing controls for the Navier-Stokes system}
In this section we give the proof of Theorem 1.1. Using similar arguments to those employed in \cite{Fer} and \cite{Fomin}, we will see that the result obtained in the previous section allows us to locally invert a nonlinear operator associated to the nonlinear system

\begin{equation} \label{34} \left\lbrace { \begin{tabular}{llll}
 $w_t- \Delta w + (w,\nabla)w+\nabla p=f+v\mathbf{ 1}_{\omega},\ \nabla \cdot w=0 $ & \textrm{in }$Q$,\\
$-z_t- \Delta z + (z,\nabla^t)w- (w,\nabla)z+\nabla q=w\mathbf{ 1}_{\mathcal O},\ \nabla \cdot z=0 $ & \textrm{in }$Q$,\\
$w=z=0$ & \textrm{on }$\Sigma$,\\
$  w_{|t=0}=0,\  z_{|t=T}=0 $ & \textrm{in }$\Omega$.\\ \end{tabular}
} \right. 
\end{equation}

%It seems natural that the insensitizing property for (\ref{34}) is a local property (in the sense of small disturbances near the equilibrium) : this explain the condition $y_0=0$ in the main theorem. Thus, the exact controllability property we want to prove for system (\ref{34}) reduce to a local exact controllability one.\\
We will use the following form of Lyusternik theorem (see \cite{Fomin}) which is in fact an inverse mapping theorem:
\begin{theoreme}
 Let $\mathcal E$ and $\mathcal G$ be two Banach spaces and let $\mathcal A:\mathcal E \mapsto \mathcal G$ satisfies $\mathcal{A}\in C^1(\mathcal{E};\mathcal{G})$. Assume that $e_0\in \mathcal E$, $\mathcal{A}(e_0)=h_0$ and $\mathcal{A}'(e_0):\mathcal E\mapsto \mathcal G$ is surjective. Then there exists $\delta>0$ such that, 
for every $h\in \mathcal G$ satisfying $\|h-h_0\|_{\mathcal G}<\delta$, there exists a solution of the equation
$$\mathcal{A}(e)=h,\quad e\in \mathcal E.$$
\end{theoreme}
 We will be led to use this theorem with the space $\mathcal{E}=\mathcal E^{s,\lambda}$, with fixed $s$ and $\lambda$ like in Theorem $3.1$ (so Lemma $4.1$ holds), 
$$\mathcal G=\mathcal G_1\times \mathcal G_2=L^2(e^{{5\over2}s\beta^{\star}} (\gamma^{\star})^{-3/2};L^2(\Omega)^N)\times L^2(e^{2s\beta^{\star}}(\gamma^{\star})^{-3/2};L^2(\Omega)^N)$$
and the operator
\begin{equation}\label{ope}\mathcal{A}(w,z,p,q,v)=(w_t- \Delta w + (w,\nabla)w+\nabla p-v\mathbf{ 1}_{\omega},-z_t- \Delta z + (z,\nabla^t)w- (w,\nabla)z+\nabla q-w\mathbf{ 1}_{\mathcal O}),\\ 
\quad \forall (w,z,p,q,v)\in\mathcal E .
\end{equation}
Since all the terms arising in the definition of $\mathcal{A}$ are linear, except for $(w,\nabla)w$ and $(z,\nabla^t)w- (w,\nabla)z$ (which are in fact bilinear), we only have to check that the terms $(w,\nabla)w$ and $(z,\nabla^t)w- (w,\nabla)z$ are well-defined and depend continuously on the data.
\begin{prop}\label{prop}
 $\mathcal{A}\in C^1(\mathcal{E};\mathcal{G}) $.
\end{prop}
$\textbf {Proof\ of\ Proposition\ \ref{prop}.}$
Let $((w^{\star},p^{\star}),(z^{\star},q^{\star}))=(e^{\frac{3}{2}s\beta^{\star}}(\widehat\gamma)^{-15/2}(w,p),e^{\frac{1}{2}s\beta^{\star}}(\widehat\gamma)^{7}(z,q))$. Then $(w^{\star},z^{\star},p^{\star},q^{\star})$ solves

\begin{equation} \label{eq33} \left\lbrace { \begin{tabular}{llll}
 $w^{\star}_t- \Delta w^{\star} + \nabla p^{\star}=f^{\star}_1+v^{\star}\mathbf{1}_{\omega}+(e^{\frac{3}{2}s\beta^{\star}}(\widehat\gamma)^{-15/2})_tw,\ \nabla \cdot w^{\star}=0 $ & \textrm{in }$Q$,\\
$-z^{\star}_t- \Delta z^{\star}+\nabla q^{\star}=f_2^{\star}+w^{\star\star}\mathbf{ 1}_{\mathcal O}-(e^{\frac{1}{2}s\beta^{\star}}(\widehat\gamma)^{7})_tz,\ \nabla \cdot z^{\star}=0 $ & \textrm{in }$Q$,\\
$w^{\star}=z^{\star}=0$ & \textrm{on }$\Sigma$,\\
$  w^{\star}_{|t=0}=z^{\star}_{|t=T}=0$ & \textrm{in }$\Omega$,\\ \end{tabular}
} \right. 
\end{equation} 
where
$$f^{\star}_1=e^{\frac{3}{2}s\beta^{\star}}(\widehat\gamma)^{-15/2}f_1,\ f^{\star}_2=e^{\frac{1}{2}s\beta^{\star}}(\widehat\gamma)^{7}f_2,\  v^{\star}=e^{\frac{3}{2}s\beta^{\star}}(\widehat\gamma)^{-15/2}v\quad \text{and}\quad w^{\star\star}=e^{\frac{1}{2}s\beta^{\star}}(\widehat\gamma)^{7}w.$$
First we look to the equation satisfied by $w^{\star}$.
We prove that the right-hand side of the first equation in (\ref{eq33}) is in $L^2(Q)^N$.
Indeed, by the definition of $\beta$, $\beta^{\star}$, $\widehat\gamma$ and $\gamma^{\star}$ we have 
\begin{itemize}
\item
$|v^{\star}\mathbf{ 1}_{\omega}|=e^{\frac{3}{2}s\beta^{\star}}(\widehat\gamma)^{-15/2}|v\mathbf{ 1}_{\omega}|\leq C(s,\lambda)e^{2s\beta^{\star}-{1\over2}s\beta}\gamma^{-15/2}|v|\mathbf{ 1}_{\omega}\in L^2(Q)^N.$
\item 
$|f_1^{\star}|=e^{\frac{3}{2}s\beta^{\star}}(\widehat\gamma)^{-15/2}|f_1|\leq C(s,\lambda)e^{{5\over 2}s\beta}(\gamma^{\star})^{-3/2}|f_1|\in L^2(Q)^N.$
\item
$|(e^{\frac{3}{2}s\beta^{\star}}(\widehat\gamma)^{-15/2})_tw|\leq CTse^{\frac{3}{2}s\beta^{\star}}(\widehat\gamma)^{-63/10}|w|\leq C(s,\lambda,T)e^{s\beta+s\beta^{\star}}(\gamma^{\star})^{-5/2}|w|\in  L^2(Q)^N$.\\

Here, we have used the fact that $e^{s\beta^{\star}}\leq C_{\epsilon} e^{s(1+\epsilon)\widehat \beta} $ for all $\epsilon>0$ and some $ C_{\epsilon}(s,\lambda)>0$.
\end{itemize}
Then, we can apply regularity results for the Stokes system (see, for instance, \cite{Temam}), hence 
\begin{equation}
 w^{\star}\in L^2((0,T);H^2(\Omega)^N)\cap L^{\infty}((0,T);H^1(\Omega)^N)\cap H^1((0,T);L^2(\Omega)^N)
\end{equation}
and depends continuously on the right-hand side of the first equation in (\ref{eq33}).
Then, if $(w,z,p,q,v)\in \mathcal E$, we have
 \begin{equation}
 e^{\frac{3}{2}s\beta^{\star}}(\widehat\gamma)^{-15/2} \nabla w\in L^{\infty}((0,T);L^2(\Omega)^{N\times N})
 \end{equation}
and 
\begin{equation}
  e^{\frac{3}{2}s\beta^{\star}} (\widehat\gamma)^{-15/2} w\in L^{2}((0,T);H^2(\Omega)^N)\subset L^{2}((0,T);L^{\infty}(\Omega)^N),
\end{equation}
thanks to the Sobolev embedding theorem. Consequently we have 
 \begin{equation}
e^{{5\over2}s\beta^{\star}}(\gamma^{\star})^{-3/2}(w,\nabla)w\leq e^{3s\beta^{\star}}(\widehat\gamma)^{-15}(w,\nabla)w\in L^2(Q)^N
 \end{equation}
and is bilinear continuous from $\mathcal E\times \mathcal E$ to $\mathcal{G}_1$.\\
Now we turn to the equation satisfied by $z^{\star}$.
%Consider $(z^{\star},q^{\star})=e^{s\beta^{\star}}(\widehat\gamma)^{d}(z,q)$. Then, $(z^{\star},q^{\star})$ satisfies (\ref{eq33}) with $d$ instead of $c$
%\begin{equation} \label{eq35} \left\lbrace { \begin{tabular}{lll}
%$-z^{\star}_t- \Delta z^{\star}+\nabla q^{\star}=w^{\star}\mathbf{ 1}_{\mathcal O}+f_2^{\star}-(e^{3s\widehat\beta-{3\over 2}s\beta^{\star}})_tz,\ \nabla \cdot z^{\star}=0 $ & \textrm{in }$Q$,\\
%$z^{\star}=0$ & \textrm{on }$\Sigma$,\\
%$z^{\star}_{|t=T}=0 $ & \textrm{in }$\Omega$.\\ \end{tabular}
%} \right. 
%\end{equation}
\begin{itemize}
\item
$|f_2^{\star}|=e^{\frac{1}{2}s\beta^{\star}}(\widehat\gamma)^{7}|f_2|\leq C(s,\lambda)e^{2s\beta^{\star}}(\gamma^{\star})^{-3/2}|f_2|\in L^2(Q)^N.$
\item
$|w^{\star\star}\mathbf{1}_{\mathcal O}|=e^{\frac{1}{2}s\beta^{\star}}(\widehat\gamma)^{7}|w|\mathbf{1}_{\mathcal O}\leq C(s,\lambda)e^{s\beta+s\beta^{\star}}(\gamma^{\star})^{-5/2}|w|\in  L^2(Q)^N.$\\

\item
$|(e^{\frac{1}{2}s\beta^{\star}}(\widehat\gamma)^{7})_tz|\leq CTse^{\frac{1}{2}s\beta^{\star}}(\widehat\gamma)^{41/5}|z|\leq C(s,\lambda,T)e^{s\beta^{\star}}|z|\in  L^2(Q)^N.$
\end{itemize}
Again, we have used the fact that $e^{s\beta^{\star}}\leq C_{\epsilon} e^{s(1+\epsilon)\widehat \beta} $ for all $\epsilon>0$ and some $ C_{\epsilon}(s,\lambda)>0$. 
 We deduce that 
\begin{equation}
z^{\star}\in L^2((0,T);H^2(\Omega)^N)\cap L^{\infty}((0,T);H^1(\Omega)^N)\cap H^1((0,T);L^2(\Omega)^N)
 \end{equation}
 and depends continuously on the right-hand side of the second equation in (\ref{eq33}).
Then, if $(w,z,p,q,v)\in \mathcal E$, we have
\begin{equation}
e^{\frac{1}{2}s\beta^{\star}}(\widehat\gamma)^{7}\nabla z\in L^{\infty}((0,T);L^2(\Omega)^{N\times N})\quad \text{and}\quad e^{\frac{1}{2}s\beta^{\star}}(\widehat\gamma)^{7} z\in L^{2}((0,T);L^{\infty}(\Omega)^{N}).
 \end{equation}
Therefore,
\begin{equation}
e^{2s\beta^{\star}}(\gamma^{\star})^{-3/2}(w,\nabla)z\in L^2(Q)^N\quad \text{and}\quad e^{2s\beta^{\star}}(\gamma^{\star})^{-3/2}(z,\nabla^t)w \in L^2(Q)^N,
 \end{equation}
since $$
(\gamma^\star)^{-3/2} \leq (\widehat\gamma)^{-1/2}.
$$

%and
%\begin{equation}
%e^{2s\beta^{\star}}(\gamma^{\star})^{-3/2}(z,\nabla^t)w\leq e^{2s\beta^{\star}}(\gamma^{\star})^{-3/2}(z,\nabla^t)w\in L^2(Q)^N.
%\end{equation}
Taking into account the continuous dependence with respect to the data, we have that these terms above are continuous from $\mathcal E\times \mathcal E$ to $\mathcal G_2$.\\
This end the proof of Proposition 5.1.\\

Finally, we can  apply  Theorem 5.1 for $e_0=0\in \mathbb {R}^{3N+2}$ and $h_0=0\in  \mathbb R^{2N}$. 
From the result obtained in Section 4, we deduce that $\mathcal A'(0,0): \mathcal E\mapsto \mathcal G$, which is given by
\begin{equation}\label{dif}
\mathcal A'(0,0)(w,z,p,q,v)=(w_t- \Delta w +\nabla p-v\mathbf{ 1}_{\omega},-z_t- \Delta z +\nabla q-w\mathbf{ 1}_{\mathcal O})\quad \forall (w,z,p,q,v)\in\mathcal E, 
\end{equation}
is surjective, that is to say $\text{Im}( \mathcal A'(0,0))=\mathcal G.$
As a conclusion, since $y_0=0$, we have find a control $v\in L^2(\omega\times (0,T))^N$ such that the associated solution to (\ref{34}) satisfies $z_{|t=0}=0$.

\end{document}